\documentclass[a4paper,reqno]{amsart}

\usepackage[algoruled,vlined,linesnumbered]{algorithm2e}
\usepackage[capitalize]{cleveref}
\usepackage[utf8]{inputenc}
\usepackage{times,eulervm}
\usepackage{enumerate}
\usepackage{booktabs}
\usepackage{amssymb}
\usepackage{siunitx}
\let\qty\SI
\usepackage{calc}
\usepackage{url}

\usepackage{ao-abbr}
\usepackage[final]{ao-ximg}
\usepackage{ao-plots}
\usepackage{ao-matrix}
\usepackage{ao-math-std}
\usepackage{ao-math-fields}
\usepackage{ao-math-symbols}
\usepackage{ao-code-names}
\usepackage{luh-colors}

\newcommand\emblack{\textcolor{black}}
\newcommand\new\emblack

\newcommand\SO[1]{\ensuremath{SO(#1)}}
\newcommand\so[1]{\ensuremath{\mathfrak{so}(#1)}}
\newcommand\setU{\mathcal U}
\newcommand\bd{\boundary}
\newcommand\qfor{\quad\forall}
\newcommand\proj[2]{#2 - \sprod{#1}{#2} #1}
\newcommand\enorm{\norm[2]}
\newcommand\mnorm{\norm[\max]}
\newcommand\orterr{e_\perp}
\newcommand\epsO[1]{\ensuremath{\eps_{O(#1)}}}
\newcommand\explog[1][]{\ensuremath{\exp_{#1}\circ\log_{#1}}\xspace}
\newcommand\explift[1][]{\ensuremath{\exp_{#1}\circ\lift_{#1}}\xspace}
\newcommand\Qty{\qty[output-exponent-marker={}]}

\DeclareMathOperator\lift{lift}
\DeclareMathOperator\rank{rank}

\DeclareMathOperator\Skew{skew}
\DeclareMathOperator\Tr{Tr}

\begin{document}

\title[PGA in Director-Based Dynamics]
{Principal Geodesic Analysis in Director-Based Dynamics\\
  of Hybrid Mechanical Systems}

\author[C. G. Gebhardt]{Cristian G. Gebhardt}
\address{Cristian G. Gebhardt\\
  University of Bergen\\
  Geophysical Institute and Bergen Offshore Wind Centre (BOW)\\
  Allégaten 70\\5007 Bergen\\Norway}
\email{cristian.gebhardt@uib.no}
\urladdr{uib.no/en/persons/Cristian.Guillermo.Gebhardt}

\author[J. Schubert]{Jenny Schubert}
\address{Jenny Schubert\\
  Leibniz Universität Hannover\\Institute of Applied Mathematics\\
  Wel\-fen\-gar\-ten 1\\30167 Hannover\\Germany}
\email{schubert@ifam.uni-hannover.de}
\urladdr{ifam.uni-hannover.de/schubert}

\author[M. C. Steinbach]{Marc C. Steinbach}
\address{Marc C. Steinbach\\
  Leibniz Universität Hannover\\Institute of Applied Mathematics\\
  Welfengarten 1\\30167 Hannover\\Germany}
\email{mcs@ifam.uni-hannover.de}
\urladdr{ifam.uni-hannover.de/mcs}

\begin{abstract}
  In this article,
  we present new computational realizations of principal geodesic analysis
  for the unit sphere $S^2$ and the special orthogonal group \SO3.
  In particular, we address the construction of long-time smooth lifts
  across branches of the respective logarithm maps.
  To this end, we pay special attention to certain critical numerical aspects
  such as singularities and their consequences on the numerical precision.
  Moreover, we apply principal geodesic analysis to investigate the behavior
  of several mechanical systems that are very rich in dynamics.
  The examples chosen are computationally modeled by employing a
  director-based formulation for rigid and flexible mechanical systems.
  Such a formulation allows to investigate our algortihms in a direct manner
  while avoiding the introduction of additional sources of error
  that are unrelated to principal geodesic analysis.
  Finally, we test our numerical machinery with the examples and,
  at the same time, we gain deeper insight into their dynamical beahvior.
\end{abstract}



\date\today

\maketitle

\section{Introduction}
\label{sec:intro}

The motion of rigid and flexible mechanical systems is inherently nonlinear
and very rich in geometric structure.
Although the mathematical theory behind the description of motion is well known,
understanding this geometric structure is by no means trivial.
In general, motion can be understood abstractly as a path on some manifold,
regardless of whether it is obtained as the solution of a set of equations
or whether it is observed or measured from a real system.
Moreover, some peculiarities or some patterns
can be as well identified or recognized from it.
Thus, having such an understanding is fundamental
to interpret the language of the motion.
However, we still need to further develop our mathematical tools
to decode such a language.
One possible strategy to accomplish this challenging task
is to employ principal geodesic analysis (PGA),
which generalizes principal component analysis (PCA)
from Euclidean spaces to nonlinear Riemannian manifolds
\cite{Jolliffe:1986,Lee:2007}.
Even though the basic ideas behind PGA are mathematically well founded
and established in the community,
exact details depend very much on the particular manifold under consideration,
and technical requirements may depend on the application context.
In general there are many ways to develop a computational realization
that is sufficiently generic, performant, robust, and accurate.
Moreover, some of these areas still remain either relatively unexplored
or not investigated.
This is the primary reason why this field of research
has attracted much attention in recent years.

As follows, we indicate some of the most important contributions to PGA
made in the last two decades.
Fletcher et al.\ \cite{Fletcher:2003} carry out statistics of shape
for medical image processing and recognize that medial descriptions
are in fact elements of a Lie group,
while Fletcher et al.\ \cite{Fletcher2004} introduce PGA
for the nonlinear statistics of shape based upon the previous findings.
In \cite{Said:2007} the exact computation of the PGA
of data on the rotation group \SO3 is presented.
Relying on PGA, Tournier et al.\ \cite{Tournier_et_al:2009}
develop a lossy compression method for human motion data.
Huckemann et al.\ \cite{Huckemann2010}
introduce a general framework for PCA on quotient spaces that result
from an isometric Lie group action on a complete Riemannian manifold.
Fotouhi and Golalizadeh \cite{Fotouhi:2012} investigate the DNA
(deoxyribonucleic acid) molecular structure by means of PGA.
Ren et al.\ \cite{Ren:2013} employ PGA to construct classifiers
for embedded non-Euclidean data.
In \cite{Sommer:2014} methods are developed
to numerically solve optimization problems over spaces of geodesics,
where PGA represents an important application of such an optimization strategy.
Salehian et al.\ \cite{Salehian:2014} present an incremental version of PGA
with applications to movement disorder classification.
Chakraborty et al.\ \cite{Chakraborty:2016}
present an efficient exact-PGA algorithm for constant curvature manifolds.
In \cite{Lazar:2017} Taylor expansions are obtained for PGA,
which lead to closed-form approximations
while revealing the impact of scale, curvature, and the distribution of data
on PGA and on the corresponding tangent space approximations.
Sittel et al.\ \cite{Sittel:2017}
specialize PGA to the torus in the context of protein dynamics.
In \cite{Heeren:2018} an invariant PGA in the space of discrete shells
is proposed and applied to model-constrained mesh editing and to
reconstruction of dense animated meshes from sparse motion capture markers.
Cazelles et al.\ \cite{Cazelles2018}
consider the analysis of datasets whose elements are random histograms and
provide a detailed comparison of PGA versus log-PCA in the Wasserstein space.
Gebhardt et al.\ \cite{Gebhardt_et_al:2019} present a PGA framework
to analyze the nonlinear dynamics of beam structures.
Curry et al.\ \cite{Curry:2019} develop a PGA specialization
based on nested sequences of totally geodesic submanifolds of symmetric spaces.
And recently, Hernandez \cite{Hernandez2021}
investigates the efficient computation of diffeomorphisms
belonging to geodesics using PGA in spaces of diffeomorphisms.

In this article, we present new computational realizations of PGA
for the unit sphere $S^2$ and the special orthogonal group \SO3.
In particular, we address the construction of long-time smooth lifts
across branches of the respective logarithm maps. This requires
paying special attention to certain critical numerical aspects
that have not been addressed extensively in the literature,
such as singularities (which arise due to the
compactness and rotational symmetry of the two manifolds)
and their consequences on the numerical precision.
Moreover, we apply PGA to investigate the dynamic behavior
of several mechanical systems:
a swinging rubber rod,
a free-oscillating cantilever beam,
a flexible triple pendulum,
a horizontal-axis wind turbine,
a tumbling cylinder,
a free-flying plate,
and a shell pendulum.
To this end, the examples are computationally modeled
by employing a director-based formulation
for rigid and flexible mechanical systems
previously developed by one of the authors
\cite{Gebhardt2017,Gebhardt2018cm,Gebhardt:2020}.
Such a director-based formulation provides a natural description of motion
in terms of product manifolds of $S^2$ or \SO3 and thus,
the application of our PGA realization is direct and straightforward.
Although more standard mechanical formulations can also be considered,
these would require further numerical pre- and post-processing of the solution.
Such additional steps would introduce, of course,
further sources of error that are not related, in any form, to PGA.
Moreover, these errors are even not segregable for our accuracy analysis
and therefore, we refrain from considering those director-free formulations.
At any rate, the choice of the examples listed before is by no means arbitrary.
On the one hand, the examples chosen exhibit very rich dyanmics
while enabling us to deeply test the numerical machinery of our PGA realization.
On the other hand, applying PGA to these examples
allows us to gain deeper insight in their dynamic behavior.
Gaining such a deeper insight is not possible with standard linear methods
such as modal analysis or PCA
due to the fact that these systematically destroy
the underlying geometric structure of the motion.
Thus, the current work provides a very robust approach to
improve the understanding of the subject.
This approach has as well the potential to be used for deriving
data-driven reduced-order models that are truly structure-preserving.
However, such a derivation of reduced-order models
falls outside the scope of the current manuscript
and will therefore be addressed in future work.

The present work is organized as follows:
In \cref{sec:dir-bas-dyn},
we brief\/ly introduce the underlying mechanical setting
relying on a director-based kinematical description.
Herein we present the \emph{canonical models}
together with their numerical treatment,
\ie, discretization in space and time.
In \cref{sec:theory}, we succinctly outline the theory behind PGA.
Then, we consider the unit sphere $S^2$
as well as the special orthogonal group \SO3,
which are two manifolds that play a fundamental role in classical mechanics.
In \cref{sec:num}, we describe the numerical implementation of PGA
for the two manifolds of interest and carry out statistics
to characterize the accuracy that can be achieved.
In \cref{sec:examples}, we carry out numerical investigations
by considering the mechanical systems mentioned above.
Finally, \cref{sec:conclusion} is devoted to conclusions,
limitations and possible future work.

\newcommand\fint{f\tsp{int}}
\newcommand\fext{f\tsp{ext}}
\newcommand\geb{\tsb{geb}}
\newcommand\ges{\tsb{ges}}
\newcommand\rb{\tsb{rb}}
\renewcommand\H{\mathcal H}
\renewcommand\K{\mathcal K}
\newcommand\setN{\mathcal N}
\newcommand\M{\mathcal M}
\newcommand\Wint{W\tsp{int}}
\newcommand\pd{\partial}

\section{Director-Based Dynamics}
\label{sec:dir-bas-dyn}

To describe the dynamic behavior of a
holonomically constrained mechanical system, we need to consider:
\emph{i)} a base manifold $Q$ (often a Euclidian space) of dimension~$s$
with a \emph{configuration map} $q\: \_\Omega \x [0, T] \to Q$
where $\Omega \subset \R^d$ for $d \in \set{1, 2, 3}$
is open and bounded with Lipschitz boundary; and,
\emph{ii)} a \emph{constraint map} $h\: Q \to \R^m$
such that $h(q(\theta, t))$ vanishes identically
for $(\theta, t) \in \_\Omega \x [0, T]$.
Then, $Q^h \define \defset{q \in Q}{h(q) = 0}$
defines the submanifold of $Q$ on which the
actual dynamics of the mechanical system is going to take place.
Note that we prefer to work on $Q$ rather than
describing the system by means of a minimal representation, \ie,
working directly on the \emph{configuration manifold} $Q^h$,
since this usually allows for a more convenient
and efficient numerical formulation.
Due to practical reasons, we adopt a variational formulation
of the dynamics based on the least action principle,
which holds for almost every $t \in [0, T]$:
\begin{equation}
  \label{eq:var-form}
  \int_\Omega \left( \sprod{\delta q}
    {\M(q) \nabla_{\!\dot{q}}(\dot{q}) + \fint(q,S) - \fext + \H(q)\tp \lambda} +
    \sprod{\delta \lambda}{h(q)} \right) d \theta = 0.
\end{equation}
For this, we also need to consider the linear image space $\Lambda$
of the Lagrange multipliers $\lambda\: \_\Omega \x [0,T] \to \Lambda$.
Then, $\delta q$ and $\delta \lambda$ are admissible variations
of $q$ and $\lambda$, respectively,
hence $\delta q(\theta) \in T_{q(\theta, t)}Q$ and
$\delta \lambda(\theta) \in T_{\lambda(\theta, t)} \Lambda \cong \Lambda$.
Moreover, $\M(q)$ is a symmetric, bounded and nonnegative mass operator.
It is the metric tensor on $Q^h$, induced by the kinetic energy functional
\begin{equation}
  \frac12 \int_\Omega \varrho \sprod{\dot{x}}{\dot{x}} \, d \theta
  =
  \frac12 \norm[\M(q)]{\dot{q}}^2,
\end{equation}
where $\varrho(\theta)$ is the mass density per unit of volume.
Moreover, $x(\theta, t) \in \R^3$ is the position of $q(\theta, t)$
in an inertial coordinate system,
and $\nabla_{\!\dot{q}}(\dot{q})$ denotes the covariant derivative
of the velocity field along itself,
which for a constant mass operator
reduces to the second time derivative of the configuration,
namely $\nabla_{\!\dot{q}}(\dot{q}) = \ddot{q}$.
The internal force covector field $\fint$ is implicitly defined
through the identity
\begin{equation*}
  \int_\Omega \sprod{\delta q}{\fint(q, S)} \, d \theta
  =
  \int_\Omega \sprod{\delta E (q)}{S} \, d \theta,
\end{equation*}
where $E(q)$ is a suitable strain measure,
\ie, invariant under affine transformations,
$S$ is the concomitant stress measure,
and $\sprod{\delta E (q)}{S}$ denotes the Euclidian scalar product
(if $(E, S)$ are represented as vectors) or the Frobenius scalar product
$\Tr(\delta E (q)\tp S)$ (if they are represented as matrices).
The pair $(E, S)$ depends on the specific mechanical model under consideration.
Later on, we are going to introduce three mechanical models,
namely the rigid body, the goemetrically exact beam,
and the geometrically exact shell.
At any rate, the former identity establishes
that the virtual work done by internal forces
is identical to the work done by internal stresses,
where the stress measure is given by
\begin{equation*}
  S = \pd_E \Wint(E)
\end{equation*}
for $\Wint$ representing an internal elastic energy functional.
Without loss of generality, the external force covector field $\fext$
is assumed to be only a function of time and bounded.
The operator $\H(q) \in L(T_qQ, \R^m)$ is the derivative of $h(q)$
and is assumed to be bounded and surjective.
At this point, it is worthwhile to point out that we are going to assume
that $\M(q)$ is positive definite on the subspace
$\ker{(\H(q))} \subseteq T_qQ$,
which together with the condition $\ker{(\M(q))} \cap \ker{(\K(q))} = \set{0}$,
for $\K(q) \in L(T_qQ, T_qQ^*)$ being the bounded and
nonnegative operator associated to the bilinear form
\begin{equation*}
  (\delta q_1, \delta q_2) \mapsto
  \pd^2_E \Wint(E(q))(\pd_q E(q) \delta q_1, \pd_q E(q) \delta q_2),
\end{equation*}
provides necessary and sufficient conditions
for well-possednes of the problem in its local form.
It should be noticed that nonnegativity, boundedness, and surjectivity
occur in the appropriate spaces, which are certainly model specific.

Before presenting the mechanical models,
we define the two manifolds
that are of special interest for our present analysis.
The first one is the Euclidian $2$-sphere
\begin{equation*}
  S^2 \define \defset{d \in \R^3}{\enorm{d} = 1},
\end{equation*}
where $d$ is normally called a \emph{director}.
The second one is the special orthogonal group
\begin{equation*}
  \SO3 \define \defset{R \in \R^{3 \x 3}}{R\tp R = I, \, \det(R) = 1},
\end{equation*}
where $R$ is the \emph{rotation tensor}.
A rotation tensor can be built from three directors as
\begin{equation*}
  R = d_1 \ox e^1 + d_2 \ox e^2 + d_3 \ox e^3,
\end{equation*}
with $\set{e^1, e^2, e^3}$ standing for the covariant basis
associated to any given orthonormal basis $\set{e_1, e_2, e_3}$
while providing that the directors' orthogonality conditions hold, namely
\begin{equation*}
  \sprod{d_1}{d_2} = \sprod{d_2}{d_3} = \sprod{d_3}{d_1} = 0.
\end{equation*}
These basic definitions are succinctly presented
with the purpose of placing the current ideas in an adequate setting.
However, further mathematical details on these two manifolds
and some interesting numerical properties
are investigated in the coming sections.

In the context of nonlinear dynamics,
there are basically three \emph{canonical models}
that attract in particular our attention.
These rely then on configuration manifols that differ from the trivial one,
namely the three-dimensional Euclidean space.
Below, we are going to introduce them brief\/ly.
For sake of brevity, we assume that
the greek indices $\alpha$ and $\beta$ can take values $\set{1, 2}$ and
the latin indices $i$ and $j$ can take values $\set{1, 2, 3}$.
It is of special convience in this section
to adopt as well the repeated-index convention.
We also use the subscript $(\fcdot)_0$
when we have adopted a specific reference,
such as $L_0$, $S_0$, or $\Omega_0$ for intital length,
initial surface, or initial volume, correspondingly,
or a reference point for the rigid body,
the center line for the beam,
or the mean surface for the shell, respectively.

The \emph{first canonical model} is the rigid body,
whose configuration manifold $Q\rb^h$ is
$\R^3 \x \SO3$ of dimension six
in $Q\rb = \R^3 \x \R^{3 \x 3}$ of dimension 12,
with configuration map $q\rb: [0, T] \to Q\rb$.
For representing \SO3, we use three directors $d_i(t)$
constrained to remain orthonormal at any time.
Then the configuration, position, and constraint maps,
with $\Omega\rb \subset \R^3$ the set of all points belonging to the rigid body
and $\theta \in \Omega\rb$, are given as
\begin{align*}
  q\rb(t)
  &= (x_0(t), d_1(t), d_2(t), d_3(t)) \in Q\rb \cong \R^{12}, \\
  x\rb(\theta, t)
  &= x_0(t) + \theta^1 d_1(t) + \theta^2 d_2(t) + \theta^3 d_3(t) \in \R^3, \\
  h\rb(x_0, d_1, d_2, d_3)
  &= (\set{\enorm{d_i}^2 - 1}_{i = 1}^3,
    \sprod{d_1}{d_2}, \sprod{d_2}{d_3}, \sprod{d_3}{d_1}) \in \R^6.
\end{align*}
As its name indicates, a rigid body is infinitely stiff and therefore,
its internal potential energy functional is the trivial one, namely
\begin{equation*}
  \Wint\rb = 0.
\end{equation*}
Clearly, the chosen setting is not minimal at all.
However, this leads to mass operators that are constant and thus,
it turns out to be very conveninent from the numerical point of view.
The resulting functions are (multivariate) polynomials of order two,
due to the presence of the orthonormality constraints.

The \emph{second canonical model} is the geometrically exact beam,
whose configuration manifold $Q\geb^h$ is again $\R^3 \times \SO3$,
with configuration map
$q\geb\: [0, L_0] \x [0, T] \to Q\geb = Q\rb$.
Any point belonging to the geometrically exact beam is represented
by the cross-sectional coordinates $(\theta^1, \theta^2) \in A_0$
and the length coordinate $\theta^3 \in [0, L_0]$ as
\begin{align*}
  q\geb(\theta^3, t)
  &= (x_0(\theta^3, t), d_1(\theta^3, t), d_2(\theta^3, t), d_3(\theta^3, t))
    \in Q\geb \cong \R^{12},
  \\
  x\geb(\theta, t)
  &= x_0(\theta^3, t) + \theta^1 d_1(\theta^3, t) + \theta^2 d_2(\theta^3, t)
    \in \R^3, \\
  h\geb(x_0, d_1, d_2, d_3) &= h\rb(x_0, d_1, d_2, d_3) \in \R^6.
\end{align*}
Even provided that a two-director formulation is possible,
we choose a three-director formulation,
which again strongly simplifies the complexity of the governing equations.
In contrast to the rigid body,
we can indeed define model-specific deformation measures
such as the axial and shear strains as
\begin{equation*}
  \Gamma^i \define \sprod{\pd_{\theta^3} x\new{_0}}{d_i}-\Gamma^i_0,
\end{equation*}
as well as the bending and torsional strains as
\begin{equation*}
  K^i \define
  \frac12 \varepsilon^i_{jk} \left(
    \sprod{\pd_{\theta^3} d_j}{d_k} - \sprod{\pd_{\theta^3} d_k}{d_j}
  \right)-K^i_0.
\end{equation*}
The internal elastic energy functional per unit of length
arises after integrating over the cross section and then,
we have the specialization
\begin{equation*}
  \int_{\Omega_0} \Wint(E) \, d \theta
  \leadsto
  \int_0^{L_0} \Wint\geb(\Gamma, K) \, d \ell.
\end{equation*}
From this, internal resultants per unit of lenght can be defined as
\begin{equation*}
  N\geb = \pd_\Gamma \Wint\geb(\Gamma, K)
  \qtextq{and}
  M\geb = \pd_K \Wint\geb(\Gamma, K),
\end{equation*}
respectively, and the components of the
internal force covector field $\fint\geb$ become
\begin{equation*}
  (\fint\geb(q, N\geb, M\geb))_a =
  (\pd_q \Gamma)^i_a (N\geb)_i + (\pd_q K)^i_a (M\geb)_i.
\end{equation*}

The \emph{third canonical model} is the geometrically exact shell,
whose configuration manifold $Q\ges^h$ is $\R^3 \x S^2$ of dimension five
in $Q\ges = \R^3 \x \R^3$ of dimension six, with configuration map
$q\ges = (x_0, d_3)\: S_0 \x [0, T] \to Q\ges$.
Any point across the geometrically exact shell is represented by
$\theta^3 \in [-\frac12 H_0, \frac12 H_0]$
(with $H_0$ standing for the thickness)
at the surface coordinates $(\theta^1, \theta^2) \in S_0$
and is described by
\begin{align*}
  q\ges(\theta^1, \theta^2, t)
  &= (x_0(\theta^1, \theta^2, t), d_3(\theta^1, \theta^2, t))
    \in Q\ges \cong \R^6,
  \\
  x\ges(\theta, t)
  &= x_0(\theta^1, \theta^2, t) + \theta^3 d_3(\theta^1, \theta^2, t)
    \in \R^3, \\
  h\ges(x_0, d_3) &= \enorm{d_3}^2 - 1 \in \R^1.
\end{align*}
Again we can define model-specific deformation measures
such as the membrane strains,
\begin{equation*}
  \epsilon^{\alpha \beta} \define
  \sprod{\pd_{\theta^\alpha} x_0}{\pd_{\theta^\beta} x_0} - \epsilon^{\alpha \beta}_0,
\end{equation*}
the bending strains,
\begin{equation*}
  \kappa^{\alpha \beta} \define
  \sprod{\pd_{\theta^\alpha} x_0}{\pd_{\theta^\beta} d_3} +
  \sprod{\pd_{\theta^\beta} x_0}{\pd_{\theta^\alpha} d_3} - \kappa^{\alpha \beta}_0,
\end{equation*}
and the transverse shear strains,
\begin{equation*}
  \gamma^\alpha \define \sprod{\pd_{\theta^\alpha} x_0}{d_3} - \gamma^\alpha_0.
\end{equation*}
To avoid confusion of the model-specific measures with those
corresponding to the geometrically exact beam, we use small fonts.
However, for the reader with experience in continuum mechanics
this would result to be a small modification
with respect to the standard notation,
where small fonts are used for the spatial description
while large fonts are used for the material description.
The internal elastic energy functional per unit of surface
arises after integrating across the thickness and then,
we have the specialization
\begin{equation*}
  \int_{\Omega_0} \Wint(E) \, d \theta
  \leadsto
  \int_{S_0} \Wint\ges(\epsilon, \kappa, \gamma) \, d S.
\end{equation*}
From this, internal resultants per unit of surface can be defined as
\begin{align*}
  N\ges &= \pd_\epsilon \Wint\ges(\epsilon, \gamma, \kappa),
  & M\ges &= \pd_\kappa \Wint\ges(\epsilon, \gamma, \kappa),
  & F\ges &= \pd_\gamma \Wint\ges(\epsilon, \gamma, \kappa),
\end{align*}
respectively, and the componets of the
internal force covector field $\fint\ges$ become
\begin{equation*}
  (\fint\ges(q, N\ges, M\ges, F\ges))_a =
  (\pd_q \epsilon)^{\alpha \beta}_a (N\ges)_{\alpha \beta} +
  (\pd_q \kappa)^{\alpha \beta}_a (M\ges)_{\alpha \beta} +
  (\pd_q \gamma)^\alpha_a (F\ges)_\alpha.
\end{equation*}

\subsection{Spatial and Temporal Discretizations}

Here we consider the spatial discretization
by means of the finite element method
and then proceed with the temporal discretization
which is designed for structure-preserving integration.

A finite element in $\R^d$ is a triplet $(K, P, \Sigma)$ where:
\emph{i)} $K$ is a non-empty closed and bounded
subset of $\R^d$ with Lipschitz boundary;
\emph{ii)} $P$ is a finite-dimensional vector space of functions on $K$; and,
\emph{iii)} $\Sigma$ is a set of linear forms $\set{\sigma_i}_{i \in \setN}$,
with $\setN \define \set{1, \dots, n}$,
that are $P$-unisolvent and known as degrees of freedom.
In addition there is a basis $\set{\phi_i}_{i \in \setN}$ for $P$
such that $\sigma_i(\phi_j) = \delta_{ij}$ for any $i, j \in \setN$,
where the functions $\phi_i$ are known as shape functions.
Let the triplet $(K, P, \Sigma)$ be a scalar-valued finite element;
if there is a set of points $\set{\theta_i}_{i \in \setN}$ in $K$
where $\sigma_i(p) = p(\theta_i)$ for any $p \in P$,
the finite element is known as Lagrange finite element and
the points are known as nodes.
Then, the property $\phi_i(\theta_j) = \delta_{ij}$
holds for any $i, j \in \setN$.
The extension to vector-valued Lagrange finite elements is carried out
simply by considering all components simultaneously.

Our spatial mesh $\mathcal{T}_h$ is then
a set of compact elements $\set{K_i}_{i \in \setN_e}$,
with $\setN_e \define \set{1, \dots, N_e}$
where $N_e$ is the number of elements, such that
\begin{equation*}
  \_\Omega = \bigcup_{i \in \setN_e} K_i, \qquad
  \operatorname{int}(K_i) \cap \operatorname{int}(K_j) = \0
  \text{ if } i \ne j, \qquad
  \bd \Omega \subset \bigcup_{i \in \setN_e} \bd K_i.
\end{equation*}
Let $V(K)$ be a Sobolev space associated to the element $K$
and $\mathcal{I}_K\: V(K) \to P$
be a local interpolation operator defined by
\begin{equation*}
  \mathcal{I}_K(v) = \sum_{i \in \setN} \sigma_i(v) \phi_i, \quad v \in V(K),
\end{equation*}
for Lagrange finite elements and $P = \Poly^k(K)$ containing the polynomials
of degree up to~$k$. Then the following estimate holds
\begin{equation*}
  \norm[L^2(K)]{v - \mathcal{I}_K(v)}
  \le
  c h_K^{\ell + 1} \norm[H^{\ell + 1}(K)]{v},
  \quad v \in V(K),
  \quad \ell \leq k,
\end{equation*}
with $c$ being a constant,
$h_K$ representing the diameter of the largest ball contained in $K$,
and all partial derivatives assumed in the weak sense.
The estimate can be extended to all of $\Omega$ as follows:
\begin{equation*}
  \norm[L^2(\Omega)]{v - \mathcal{I}_h(v)}
  \le
  c h^{\ell + 1} \norm[H^{\ell + 1}(\Omega)]{v},
  \quad v \in V,
  \quad \ell \leq k,
\end{equation*}
where $\mathcal{I}_h$ is the global interpolation operator,
$h$ is the largest $h_K$ in the mesh $\mathcal{T}_h$
and $V$ is the space in which the variational problem is formulated.
For flexible mechanical models considered in this work we need to consider
\begin{equation*}
  V \define H^1(\Omega)^s,
\end{equation*}
the continuous space equipped with $s$ coordinate fields $q_i$
in which the problem is stated, and
\begin{equation*}
  V_h \define
  \defset{v_h \in C^0(\Omega)^s \cap V}
  {v_h|_{K_i} \in \Poly^1(K_i) \, \forall i \in \setN_e},
\end{equation*}
the discrete space in which the problem is approximated.
Thus, $V_h$ is a subspace of $V$.

Particularly for our two flexible models,
the geometrically exact beam and geometrically exact shell, we have
\begin{equation*}
  q\geb(\theta^3; t) \approx q\tsb{geb,h}(\theta^3; t) =
  \sum_{i \in \setN_e} \Phi_{\text{geb},i}(\theta^3) \^q_{\text{geb},i}(t)
  \qquad (s = 12)
\end{equation*}
and
\begin{equation*}
  q\ges(\theta^1,\theta^2; t) \approx q\tsb{ges,h}(\theta^1,\theta^2; t) =
  \sum_{i \in \setN_e} \Phi_{\text{ges},i}(\theta^1, \theta^2) \^q_{\text{ges},i}(t)
  \qquad (s = 6),
\end{equation*}
where $\Phi_i$ is a matrix containing all shape functions
and $\^q_i$ is the set of nodal degrees of freedom.
The two-node beam element has twelve coordinate fields per node,
one position vector and three directors,
$\^q = (x_{0,i}, d_{1,i}, d_{2,i}, d_{3,i})_{i = 1}^2$,
whereas the four-node shell element has six coordinate fields per node,
one position vector and one director,
$\^q = (x_{0,i}, d_{3,i})_{i = 1}^4$.
To complete the description,
we enforce the constraints associated to the directors only at the nodal level.
Then we can proceed with approximations
for all the terms present in \eqref{eq:var-form} as
\begin{align*}
  \int_\Omega \sprod{\delta q_h}
  {\M(q_h) \nabla_{\!\dot{q}_h} (\dot{q}_h)} \, d \theta
  &= \Sprod{\delta \^q}
    {\int_\Omega \^\M(\^q) \^\nabla_{\!\dot{\^q}}(\dot{\^q}) \, d \theta},
  \\
  \int_\Omega \sprod{\delta q_h}{\fint(q_h) - \fext} \, d \theta
  &= \Sprod{\delta \^q}
    {\int_\Omega \bigl( \^f\tsp{int}(\^q) - \^f\tsp{ext} \bigr) d \theta},
  \\
  \int_\Omega \sprod{\delta q_h}{\H(q_h)\tp \lambda_h} \, d \theta
  &= \Sprod{\delta \^q}{\^\H(\^q)\tp \^\lambda},
  \\
  \int_\Omega \sprod{\delta \lambda_h}{h(q_h)} \, d \theta
  &= \Sprod{\delta \^\lambda}{\^h(\^q)}.
\end{align*}
It is worthwhile to mention that upon linearization,
removing of rigid body motions from $V$ and $V_h$,
and turning off the term corresponding to the time derivative
of the momentum, the convergence rate
\begin{equation*}
  \norm[L^2(\Omega)]{q - q_h} \le c h^2
\end{equation*}
can be derived.
Even provided that this rate does not necessarily hold
for the nonlinear setting,
numerical experiments suggest that the nonlinear problem
may also exhibit optimal convergence behavior.

To implicitly integrate the system's dynamics in time,
we are going to rely on a discrete version of the balance equation
at time instant $t_{n+\frac12}$.
Moreover, the time grid is uniform
and characterized by the fixed time step $\Delta t$.
Since the problem is already discretized in space,
we assume that all quantities are finite-dimensional or nodal quantities
and therefore drop $\^{(\fcdot)}$ from our notation.
To properly set up the numerical integration scheme,
we employ second-order approximations
for the postion and the velocity,
\begin{equation*}
  q_{n+\frac12} \approx \frac{q_n + q_{n+1}}{2}
  \qtextq{and}
  \dot{q}_{n+\frac12} \approx \frac{q_{n+1} - q_n}{\Delta t},
\end{equation*}
for the time derivative of the momentum,
\begin{equation*}
  \left( \M(q) \nabla_{\!\dot{q}}(\dot{q}) \right)_{n+\frac12}
  \approx
  \frac{\M(q_{n+1}) \dot{q}_{n+1} - \M(q_{n}) \dot{q}_{n}}{\Delta t} -
  \frac12 \pd_\xi \! \left.
    \sprod{\dot{q}_{n+\frac12}}{\M(\xi) \dot{q}_{n+\frac12}}
  \right|_{\xi = q_{n+\frac12}} \!,
\end{equation*}
and for the internal force covector field,
\newcommand{\contract}{\mathbin{\raisebox{\depth}{\scalebox{1}[-1]{$\lnot$}}}}
\begin{equation*}
  \left( \fint(q) \right)_{n+\frac12}
  \approx
  \frac12 \int_{-1}^{+1}
  \left( \pd_q E \right)_{n+\frac12} \contract S(q(\xi)) \, d \xi
\end{equation*}
(where ``$\contract$'' denotes the contraction of a $(m,n+1)$-tensor
with a $(n,m)$-tensor resulting in a one-form \cite{Heard2006}) with
\begin{equation*}
  q(\xi) \define \frac12 (1 - \xi) q_n + \frac12 (1 + \xi) q_{n+1}
\end{equation*}
and
\begin{equation*}
  \left( \pd_q E \right)_{n+\frac12} = \pd_q E(q_{n+\frac12}).
\end{equation*}
In addition, the constraint force is approximated as
\begin{equation*}
  (\H(q)\tp \lambda)_{n+\frac12}
  \approx
  \frac12 \int_{-1}^{+1} \H(q(\xi))\tp \lambda_{n+\frac12} \, d \xi
\end{equation*}
while the constraint equation is enforced at time instant $t_{n+1}$,
\ie $h(q_{n+1}) = 0$.

The integration scheme then proceeds as follows:
\emph{i}) $q_n$ and $\dot{q}_n$ are given and thus,
used to initially estimate $q_{n+1}$ and $\dot{q}_{n+1}$;
\emph{ii}) from the balance equation at time instant $t_{n+\frac12}$
and the constraint equation at time instant $t_{n+1}$,
$q_{n+1}$ and $\lambda_{n+\frac12}$ are to be found,
up to a tolerance $\varepsilon$,
through an iterative  procedure such as Newton's method;
and, \emph{iii}) $\dot{q}_{n+1}$ is finally calculated
with the fromula $\frac{2}{\Delta t}(q_{n+1} - q_n) - \dot{q}_n$.

This hybrid numerical strategy, beyond being second-order accurate,
identically preserves linear and angular momentum as well as total energy.
Therefore, this turns out to be good enough for our purpose.

\section{Theoretical Background of PGA}
\label{sec:theory}

PGA consists of three steps:
\emph{i)}~lifting a discrete trajectory from the manifold to a suitable
tangent space by means of \emph{logarithm maps};
\emph{ii)}~performing a PCA in that tangent space; and
\emph{iii)}~mapping the results back to the manifold
by an \emph{exponential map}.
This section provides the mathematical notation, concepts,
and facts that we need to formulate steps \emph{i)} and \emph{iii)}
for the manifolds of interest
and to understand certain difficulties caused by singularities.
The presentation is essentially rigorous
but we avoid unnecessary abstraction, technical details, and proofs.

\subsection{Riemannian Manifolds}
\label{sec:manifolds}

Before turning to the manifolds of interest,
we give a brief introduction into the general theory.
Background material can be found in any textbook on Riemannian geometry,
such as \cite{Klingenberg:1982}.
For the purpose of this paper, a \emph{(smooth) manifold} $M$
of dimension $m$ is a subset of $\R^n$
defined by a smooth surjective map $F\: \R^n \to \R^{n - m}$,
\begin{align*}
  M &\define \defset{x \in \R^n}{F(x) = 0}, &
  \rank F'(x) &= n - m \qfor x \in M.
\end{align*}
Of course, $\R^n$ is itself a manifold (set $F \define 0$).
One of the simplest nontrivial examples is the
$(n-1)$-dimensional manifold of unit vectors
(or \emph{directors}) in $\R^n$: the unit sphere
\begin{align*}
  S^{n-1} &\define \defset{d \in \R^n}{\enorm{d}^2 - 1 = \sprod{d}{d} - 1 = 0}.
\end{align*}
Every point $x$ on a manifold $M$ has an attached \emph{tangent space}
$T_xM \subset M \x \R^n$ defined as
\begin{align*}
  T_xM
  &\define
  \defset{(x, v) \in \set{x} \x \R^n}{F'(x) v = 0}
  =
  \set{x} \x \defset{v}{F'(x) v = 0}
  ,
\end{align*}
which is a \emph{linear} space isomorphic to $\R^m$.
If $I$ is an interval and $c\: I \to M$ a smooth curve,
then $\dot c(t)$ is a tangent vector in $T_{c(t)}M$ for each $t \in I$.
Moreover, the tangent space $T_xM$ depends smoothly on $x$
so that $\dot c$ is a smooth curve in the set of all tangent spaces,
called the \emph{tangent bundle} $TM$ of $M$,
which is also a manifold (of dimension $2 m$ in $\R^{2 n}$),
\begin{equation*}
  TM \define \defset{(x, v) \in \R^n \x \R^n}{(F(x), F'(x) v) = 0}.
\end{equation*}
In case of the $(n - 1)$-sphere, for instance, the tangent space at $d$ is
essentially the space $\set{d}^\perp$ of vectors that are orthogonal to $d$,
\begin{equation*}
  T_dS^{n-1}
  =
  \set{d} \x \defset{v \in \R^n}{\sprod{d}{v} = 0}
  =
  \set{d} \x \set{d}^\perp
  .
\end{equation*}
Now a \emph{Riemannian metric} on $M$ is intuitively a set of scalar products,
one on each tangent space $T_xM$, that vary smoothly with $x \in M$.
Formally it is a smooth map $g$ from $M$ into
a subset of the \emph{tensor bundle} $\mathbb S^2M$
of symmetric bilinear forms on $T_xM$,
which has a definition similar to the bundle $TM$ and is again a manifold.
Note that every scalar product on $\R^n$
is a Riemannian metric that does not depend on $x$
and that \emph{induces} such a metric on~$M$:
$g(x)(v,w) = \sprod{v}{w}$ for
$v, w \in \R^n$ or $v, w \in T_xM$.
With a Riemannian metric, $M$ becomes a \emph{Riemannian manifold}
and we can generalize the notion of straight lines in $\R^n$
to \emph{geodesics} on $M$,
which are solutions of a certain second order differential equation
that is entirely determined by the metric~$g$.
Every geodesic $c\: I \to M$ has constant speed
$s = g(c(t))(\dot c(t), \dot c(t))^{1/2}$,
and it is the unique shortest curve between
sufficiently close points $x_1 = c(t_1)$, $x_2 = c(t_2)$,
thus defining $d_M(x_1, x_2) \define s \abs{t_1 - t_2}$
as their \emph{distance} in~$M$.
Considering $S^{n-1}$ again,
we take the scalar product on $\R^n$ as Riemannian metric $g$.
Then every geodesic $c\: \R \to S^{n-1}$ has the form
$c(t) = \cos(st) d + \sin(st) e$,
where $d, e$ are orthogonal vectors in $S^{n-1}$
and $s \ge 0$ is the speed.
If $s > 0$, $c$ wraps around $S^{n-1}$ on a circular orbit.

We are now ready to address the key concepts for PGA.
To define the exponential map, we use the unique geodesic $c_v$
that has a given vector $v \in T_xM$ as initial tangent:
\begin{equation*}
  \exp_x\: \setU_x \subseteq T_xM \to M, \quad
  \exp_x(v) \define c_v(1) \qtextq{where}
  \dot c_v(0) = v.
\end{equation*}
Clearly, $\exp_x(0) = x$, and by the Picard--Lindelöf theorem
$\exp_x$ is defined on some open neighborhood $\setU_x$ of $0$.
On our manifolds of interest it is defined globally
and is actually a surjective map.
In particular, for $\exp_d\: T_dS^{n-1} \to S^{n-1}$
(with $\norm{}$ induced by $g$) we have
\begin{equation*}
  \exp_d(v) = \cos(\norm{v}) d + \sin(\norm{v}) \frac{v}{\norm{v}}
  \qfor v \in T_dS^{n-1} \without{0}.
\end{equation*}
The exponential map is always smooth,
and by the inverse function theorem
it possesses a smooth local inverse
defined on some open neigborhood $M_x$ of $x \in M$.
This inverse is called the \emph{logarithm map}, $\log_x\: M_x \to T_xM$.
Moreover, further logarithm maps may exist as local inverses of $\exp_x$
on neighborhoods of other image points $\exp_x(v) \notin M_x$.

For PGA it is crucial that $\exp_x$ is globally defined and surjective,
and that local inverses exist almost everywhere on $M$:
given a smooth curve $c\: [0, 1] \to M$ starting at $c(0) = x$,
we seek a smooth curve $y\: [0, 1] \to T_xM$
such that $exp_x(y(t)) = c(t)$ for all $t \in [0, 1]$.
To construct that \emph{lifted} curve $y$,
we need one or more logarithm maps such that
the closures of their domains of definition
cover the image of $c$.
Concrete representations of logarithm maps and lift maps
for the manifolds of interest will be given in
\cref{sec:theory:dir+rot,sec:theory:lifts}.

All the concepts above remain valid
when intersecting $M$ with an open subset $\setU$ of $\R^n$.
This still gives quite special submanifolds of $\R^n$,
and we will use it below.
More generally, all concepts extend to manifolds that are
patched together from overlapping pieces of the form $M \cap \setU$,
even without a surrounding linear space,
and even in infinite dimension.

\subsection{Directors and Rotations}
\label{sec:theory:dir+rot}

Our first manifold of interest is the set of directors in $\R^3$,
the unit 2-sphere $S^2 = \defset{d \in \R^3}{\enorm{d} = 1}$,
which has already served as an illustrative example.
Recall the properties
\begin{align}
  T_dS^2 &=
  \set{d} \x \defset{v \in \R^3}{\sprod{d}{v} = 0} \cong \set{d}^\perp, \notag \\
  \label{eq:s2:expd}
  \exp_d(v) &=
  \begin{cases}
    d, & \norm{v} = 0, \\
    \displaystyle
    \cos(\norm{v}) d + \sin(\norm{v}) \frac{v}{\norm{v}}, & \norm{v} > 0.
  \end{cases}
\end{align}
Clearly we have $\exp_d(v) = -d$
for every tangent vector of length $\norm{v} = \pi$,
and the open disk
$D_\pi \define \defset{v \in T_dS^2}{\norm{v} < \pi}$
is actually a maximal domain of injectivity of $\exp_d$.
Now the logarithm map $\log_d\: S^2 \without{-d} \to D_\pi$
has the explicit representation
\begin{align}
  \label{eq:s2:logd}
  \log_d(e) &= (d, v), &
  v &\define
  \begin{cases}
    0, & \enorm{\proj{d}{e}} = 0,\\
    \displaystyle
    \arccos \sprod{d}{e} \frac{\proj{d}{e}}{\enorm{\proj{d}{e}}}, &
    \enorm{\proj{d}{e}} > 0,
  \end{cases}
\end{align}
where $\proj{d}{e}$ is the projection of $e$ on $\set{d}^\perp$
and $\norm{v} = \arccos \sprod{d}{e} \in (0, \pi)$
is the angle between $d$ and $e$.
Further logarithm maps and the construction of a lift map
will be discussed in \cref{sec:theory:lifts}.

Our second manifold of interest is the
matrix group of \emph{rotations} in $\R^3$,
called the \emph{special orthogonal group} \SO3.
This is a (non-commutative) multiplicative group
and also a 3-dimensional manifold in the matrix space $\R^{3 \x 3}$
seen as the vector space~$\R^9$.
To understand the basic structure,
we first define the \emph{orthogonal group}
\begin{equation*}
  O(3) \define \defset{R \in \R^{3 \x 3}}{R\tp R = I}.
\end{equation*}
When viewing the matrix $R$ as a triple of mutually orthogonal directors,
$R = [d_1 \ d_2 \ d_3]$,
the symmetric matrix constraint $R\tp R = I$
is equivalently written as
\begin{align*}
  \enorm{d_1}^2 = \enorm{d_2}^2 = \enorm{d_3}^2 &= 1, &
  \sprod{d_1}{d_2} = \sprod{d_1}{d_3} = \sprod{d_2}{d_3} = 0.
\end{align*}
Thus $O(3)$ is a 3-dimensional submanifold
in the 9-dimensional matrix space $\R^{3 \x3} \cong \R^9$.
The determinant of an orthogonal matrix is either $+1$ or $-1$,
and hence $O(3)$ consists of two disconnected parts.
Now let $\setU_+ \define \defset{R \in \R^{3 \x 3}}{\det(R) > 0}$
(which is open) and define
\begin{equation*}
  \SO3 \define O(3) \cap \setU_+ =
  \defset{R \in \R^{3 \x 3}}{R\tp R = I, \det(R) = 1}.
\end{equation*}
This is still a manifold with a group structure, called a \emph{Lie group}.
On Lie groups it suffices to consider as a ``universal'' tangent space
the one at the neutral element, the \emph{Lie algebra},
since every tangent space has a natural representation in terms of it.
In our case, that representation at the identity matrix $I$ reads
\begin{equation*}
  T_R\SO3
  =
  \defset{R N}{N \in A(3)}
  =
  R T_I\SO3
  =
  R \so3
  ,
\end{equation*}
where $A(3)$ denotes the space of
\emph{antisymmetric} (or \emph{skew-symmetric}) matrices,
$A(3) \define \defset{N \in \R^{3 \x 3}}{N\tp = -N}$,
and \so3 is the standard notation for the Lie algebra of \SO3.
We also use the more succinct notation
$R N$ for $(R, R N) \in T_R\SO3$.
As Riemannian metric on $\R^{3 \x 3}$ and hence on \SO3
we choose a scaled Frobenius scalar product
(so that $\norm{N}$ becomes the rotation angle of $\exp(N)$, see below),
\begin{align*}
  \sprod{V}{W} &\define \frac{\sprod[F]{V}{W}}{2} = \frac{\Tr(V\tp W)}{2}, &
  g(R)(R N_1, R N_2) &= \sprod{R N_1}{R N_2} = \sprod{N_1}{N_2}.
\end{align*}
Thus, $T_R\SO3$ is a subspace of the 8-dimensional
orthogonal complement of $\operatorname{span}\set{R}$,
\begin{equation*}
  \set{R}^\perp =
  \defset{R V \in \R^{3 \x 3}}{0 = \sprod{R}{R V} = \tfrac12 \Tr(V)}.
\end{equation*}
Now the exponential map $\exp\: \so3 \to \SO3$
is the standard matrix exponential,
and due to the identity $N^3 = -\norm{N}^2 N$
(with $\norm{}$ induced by $g$)
it has the explicit representation
\begin{equation}
  \label{eq:so3:exp}
  \exp(N)
  =
  \sum_{k = 0}^\oo \frac1{k!} N^k
  =
  \begin{cases}
    I, & \norm{N} = 0, \\
    \displaystyle
    I + \frac{\sin \norm{N}}{\norm{N}} N +
    \frac{1 - \cos \norm{N}}{\norm{N}^2} N^2, & \norm{N} > 0.
  \end{cases}
\end{equation}
It is globally defined and surjective.
With $\theta \define \norm{N}$ and $U \define N / \theta$ for $\theta > 0$,
the last formula explains the so-called \emph{axis-angle representation}
$N = \theta U$ of the rotation $\exp(N)$,
\begin{equation}
  \label{eq:so3:exp-aa}
  \exp(N) = \exp(\theta U) = I + (\sin \theta) U + (1 - \cos \theta) U^2.
\end{equation}
To understand this fundamental geometric interpretation,
consider a rotation axis given by a director $u \in S^2$
and a rotation angle $\theta \in \R$.
Then the rotation about axis $u$ with angle $\theta$
is given by the \emph{rotation vector} (or \emph{Euler vector})
$n = \theta u$ as $\exp(N)$ where
\begin{equation*}
  N = \theta U = \theta \Skew(u) = \Skew(n) \define
  \mat[3]{0 & -n_3 & n_2 \\ n_3 & 0 & -n_1 \\ -n_2 & n_1 & 0}.
\end{equation*}
Note that our definition of $g$
makes the isomorphism $\Skew\: \R^3 \cong A(3) \cong \so3$
\emph{isometric}, \ie, $\norm{N} = \enorm{n} = \theta$.
We identify $U$ with the unit axis $u$ by virtue of this isometry
and note that $U v = u \x v$ for each $v \in \R^3$,
where $u \x v$ denotes the cross product.

Using the notation $R(U, \theta) = \exp(\theta U)$,
the identity \eqref{eq:so3:exp-aa} readily implies several
geometrically obvious symmetries of \SO3
from which much of the subsequent theory derives:
\begin{align*}
  R(-U, -\theta) &= R(U, \theta), \\
  R(U, \theta + 2 k \pi ) &= R(U, \theta) \qfor k \in \Z,
\end{align*}
and the direct implications
\begin{align*}
  R(-U, \pi) = R(U, -\pi) &= R(U, \pi), \\
  R(-U, \pi - \theta) = R(U, \theta - \pi) &= R(U, \pi + \theta).
\end{align*}
The geodesics $c\: \R \to \SO3$ are precisely
the curves $c(t) \define R \exp(t N)$
with arbitrary $N \in \so3$ and  $c(0) = R$.
The velocitiy is $\dot c(t) = c(t) N$
so that $\dot c(0) = R N$.
All geodesics are closed curves with period $2 \pi / \norm{N}$,
except for the constant ones (with $N = 0$).

A maximal domain of injectivity of $\exp$ is the open ball
$B_\pi \define \defset{N \in \so3}{\norm{N} < \pi}$
whereas every pair of antipodal points $\pm \pi U \in \bd B_\pi$
maps to the same rotation $\exp(\pi U)$.
Hence we have a unique local inverse, the logarithm map
$\log\: \SO3 \setminus \exp(\bd B_\pi) \to B_\pi$,
which is given by the explicit representation
\begin{equation}
  \label{eq:so3:log}
  \log(R)
  =
  \begin{cases}
    0, & R = I, \\
    \displaystyle
    \frac{\theta}{\sin \theta} \frac{R - R\tp}{2},
    \quad
    \theta \define \arccos \frac{\Tr(R) - 1}{2}, & R \ne I.
  \end{cases}
\end{equation}
Here $\theta \in [0, \pi)$ is again the rotation angle.
We finally note that $\exp_R(R N) = R \exp(N)$ and
$\log_R(R S) = R \log(S)$ for all $R, S \in \SO3$ and $N \in \so3$.
Further logarithm maps and the construction of a lift map
will again be discussed in \cref{sec:theory:lifts}.

\subsection{Singularities and Lifts}
\label{sec:theory:lifts}

Here we take a deeper look at the exponential maps
and associated logarithm maps on $S^2$ and \SO3.
The issues to be discussed are caused
by the \emph{rotational symmetries} and
by the \emph{compactness} that both manifolds exhibit:
the respective exponential maps ``wrap their
tangent space around the manifold'' infinitely often.
Therefore the inverses of the exponential maps are not unique
(in fact they have infinitely many possible values at certain singularities)
and hence we may need several logarithm maps
to lift a curve by ``unwrapping'' the tangent space.

\begin{figure}[tp]
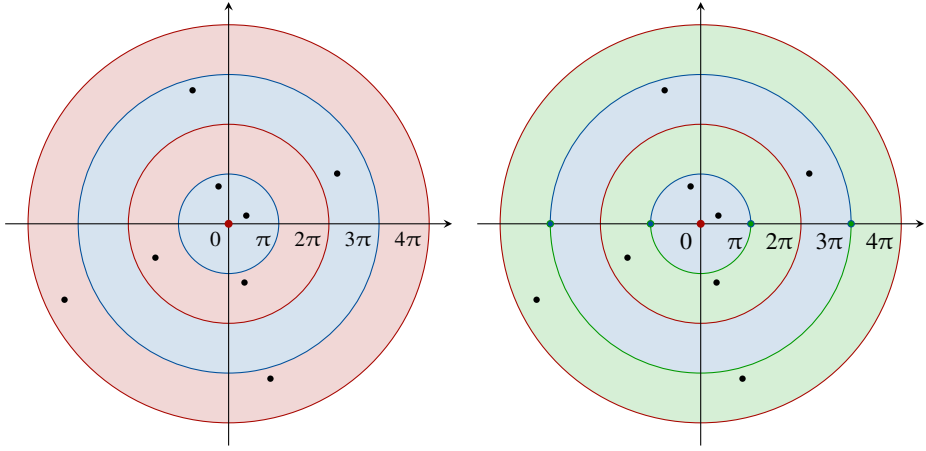

  \centering
  \AOtikz{annuli}\hfil
  \AOtikz{shells}
  \caption{Tangent space $T_dS^2$ (left) and a 2d subspace of \so3 (right).
    Light colored rings show maximal injectivity sets of $\exp_d$ \resp $\exp$.
    Red circles show ``even'' boundaries that are mapped to $d$ \resp $I$;
    blue circles (left) and antipodal blue-green half circles (right)
    show ``odd'' boundaries that are mapped to $-d$ (left)
    or $\exp(\bd B_\pi)$ (right).
    Two sets of black dots on radial lines show
    $\log_0, \dots, \log_3$ of two elements of $S^2$ \resp \SO3.}
  \label{fig:shells}
\end{figure}

Turning first to the 2-sphere $S^2$,
we fix an element $d$ and denote by $D_r$
the open circular disk of radius $r$ in $T_dS^2$.
Next we define annuli $A_0 \define D_\pi \!\without 0$ and
$A_k \define D_{(k+1) \pi} \!\setminus \_D_{k \pi}$ for $k > 0$,
where $\_D_r$ denotes the closure of $D_r$.
Then every tangent vector $v \in A_k$
has a unique axis-angle representation $v = \theta u$ with $\norm{u} = 1$
and $\theta = \norm{v} \in (k \pi, (k+1) \pi)$,
and hence we can define smooth bijections $\Phi_{k\ell}\: A_k \to A_\ell$
such that $\exp_d(\Phi_{k\ell}(v)) = \exp_d(v)$:
\begin{equation*}
  \Phi_{k\ell}(v) = \Phi_{k\ell}(\theta u) \define
  \begin{cases}
    [(\ell - k    ) \pi + \theta]   u , & \ell - k \text{ even}, \\
    [(\ell + k + 1) \pi - \theta] (-u), & \ell - k \text{ odd}.
  \end{cases}
\end{equation*}
In \cref{sec:theory:dir+rot} we have already seen that $\exp_d$
maps $D_\pi$ bijectively to $S^2 \without{-d}$
and the entire boundary $\bd D_\pi$ to $-d$.
More generally, $exp_d$ maps every annulus $A_k$
bijectively to $S^2 \without{-d, +d}$,
every ``even'' boundary $\bd D_{2 k \pi}$ to $+d$,
and every ``odd'' boundary $\bd D_{(2 k + 1) \pi}$ to $-d$.
This yields infinitely many distinct \emph{branches} of the logarithm map
(see \cref{fig:shells}),
\begin{align*}
  \log_{d,k}\: S^2 \without{\pm d} &\to A_k, &
  \log_{d,k}(e) \define \Phi_{0k}(\log_d(e)),
\end{align*}
where $\log_{d,0}$ agrees with the \emph{principal branch} $\log_d$
except that $\log_{d,0}(d)$ is not defined.
Clearly, the closures $\_A_k$ cover the tangent space,
$T_dS^2 = \bigcup_{k \in \N} \_A_k$,
and hence it should be possible to lift
smooth curves from $S^2$ to $T_dS^2$
by switching from $\log_{d,k}$ to $\log_{d,k\pm1}$
whenever the lifted curve $y(t)$ crosses $\bd A_k$.
However, $\log_{d,k}$ has singularities at $\pm d$
where the projection $\proj{d}{e}$
in \eqref{eq:s2:logd} vanishes and
the image point in $\bd A_k$ is not unique.
These singularities can prevent the construction of a smooth lift.
In fact, a curve that enters $-d$ at zero speed
along a geodesic circle and exits along another one
will have a discontinuous lift,
and a curve that stays at $-d$ for a positive time
will have an undefined lift.
The same happens at $+d$ unless the lifted point in $T_dS^2$ is the origin.
As a sufficient condition for a unique smooth lift
we require that $y(t)$ crosses $\bd A_k$ at nonzero speed,
\ie, $\dot c(t) \ne 0$ whenever $c(t) = \pm d$.
In that case, $\dot y(t)$ is necessarily orthogonal to $\bd A_k$.

On \SO3 we face an essentially similar
though slightly more complicated situation.
By $B_r$ we denote the open ball of radius $r$ in \so3.
Then we define the open spherical shells
$S_0 \define B_\pi \!\without 0$ and
$S_k \define B_{(k+1)\pi} \!\setminus \_B_{k\pi}$ for $k > 0$.
Using the unique angle-axis representation $N = \theta U$
of a matrix $N \in S_k$, where $\theta \in (k \pi, (k+1) \pi)$,
we have again smooth bijections $\Phi_{k\ell}\: S_k \to S_\ell$
such that $\exp(\Phi_{k\ell}(\theta U)) = \exp(\theta U)$:
\begin{equation*}
  \Phi_{k\ell}(N) = \Phi_{k\ell}(\theta U) \define
  \begin{cases}
    [(\ell - k    ) \pi + \theta]   U , & \ell - k \text{ even}, \\
    [(\ell + k + 1) \pi - \theta] (-U), & \ell - k \text{ odd}.
  \end{cases}
\end{equation*}
Now $\exp$ maps $\bd B_{k\pi}$ entirely to $I$ for even $k$ and maps
antipodal points on $\bd B_{k\pi}$ to the same image point for odd $k$.
Finally $S_k$ is bijectively mapped to $\SO3 \setminus \exp(\bd S_k)$,
and we obtain again a sequence of distinct branches of the logarithm map,
\begin{equation*}
  \log_k\: \SO3 \setminus \exp(\bd S_k) \to S_k, \qquad
  \log_k(R) \define \Phi_{0k}(\log(R)).
\end{equation*}
Again $\log_0$ is undefined at $R = I$ but otherwise
identical to the \emph{principal branch} $\log$.
For an illustration in $\R^3 \cong \so3$ see \cref{fig:shells}.
Note that each image set $\exp(\bd S_k)$
contains the identity~$I$ as image of the entire ``even'' boundary.
In contrast, the image of each ``odd'' boundary
is the 2-dimensional submanifold $\Pi \define \exp(\bd B_\pi)$ in \SO3:
the geodesic sphere of all rotations at distance $\pi$ to $I$.
Moreover, it is immediate from \eqref{eq:so3:exp-aa}
that $\Pi$ is also a sphere \new{at Euclidian distance} 2
\new{from} $I$ in $\R^{3 \x 3}$:
\begin{equation*}
  \norm{\exp(\pi U) - I}^2
  =
  \norm{2 U^2}^2
  =
  2 \Tr(U\tp U\tp U U)
  =
  2 \Tr(U\tp U)
  =
  4 \norm{U}^2
  =
  4
  .
\end{equation*}
Thus, both $S^2$ and \SO3 have
geodesic diameter $\pi$ and Euclidian diameter 2
(as respective subsets of $\R^3$ and $\R^{3 \x 3}$),
but every point $d$ in $S^2$ has
the single antipodal point $-d$
whereas every matrix $R$ in \SO3 has
a 2-sphere of most distant matrices: \new{the set} $R \Pi$
\new{(which can be shown to be centered
  at $-\frac13 R$ with radius $2 / \sqrt3$)}.
In fact, the antipodal matrix $-R$
is always an element of $O(3)$ outside \SO3,
at Euclidian distance $\sqrt6$.

Finally, as before, the closures of the domains of the logarithm maps
cover the entire tangent space,
$\so3 = \bigcup_{k \in \N} \_S_k$,
but we have singularities at the boundaries
where the denominator $\sin \theta$ in \eqref{eq:so3:log} vanishes
and the image in $\bd S_k$ is not unique.
To obtain unique and smooth lifts,
we again require $\dot c(t) \ne 0$
whenever $c(t) \in \exp(\bd S_k) = \set{I} \cup \Pi$.

Note that, due to the rotational symmetry of $S^2$,
we might also lift a curve $c\: I \to S^2$ into \so3
by setting $c(t) \enifed R(t) c(0)$.
This defines a unique curve $R\: I \to \SO3$
that gives the desired lift $N\: I \to \so3$.
Then $c$ is recovered as $c(t) = \exp(N(t)) c(0)$.
We do not explore this relation further;
the underlying mathematical property is that
the symmetric space $S^2$ becomes a homogeneous space
under the group action of $\SO3$.

We conclude this section by noting that the following
orthogonal projections are defined for $S^2$ and \SO3
as respective submanifolds of $\R^3$ and $R^{3 \x 3}$:
given $x \in \R^3 \without{0}$, the orthogonal projection is
\begin{equation}
  \label{eq:s2:proj}
  \frac{x}{\enorm{x}} \in S^2.
\end{equation}
Given $B \in \R^{3 \x 3}$ with $\det(B) > 0$ and with
singular value decomposition $B = U \Sigma V\tp$,
the orthogonal projection is
\begin{equation}
  \label{eq:so3:proj}
  U V\tp \in \SO3.
\end{equation}

\section{Numerical Implementation of PGA}
\label{sec:num}

In this section we discuss the numerial implementations
of the exponential map and of the principal branch of the logarithm map
for our manifolds of interest:
\begin{align*}
  \exp_d\: &T_dS^2 \to S^2, &
  \exp\: &\so3 \to \SO3, \\
  \log_d\: &S^2 \without{-d} \to D_\pi \subset T_dS^2, &
  \log\: &\SO3 \setminus \Pi \to B_\pi \subset \so3.
\end{align*}
For the logarithms we have to take special care
in the vicinity of their singularities.
Later on we construct the corresponding lifts
of discrete snapshots $x_i = c(t_i)$
of a smooth curve $c\: [t_0, t_n] \to M$ for $M = S^2$ and $M = \SO3$.
We assume that $c$ itself is not given
and define the discrete lift $y_0, \dots, y_n$ as follows.
Set $y_0 \define 0$ in $T_{x_0}M$.
For $i = 1, \dots, n$, lift the geodesic segment
from $x_{i-1}$ to $x_i$ into $T_{x_0}M$
so that the lift starts at $y_{i-1}$,
and define its endpoint as $y_i$.
We call this step the \emph{lift map},
denoted as $\lift_{x_0, y_{i-1}}(x_i)$.
On $S^2$ and \SO3, the geodesic segment from $x_{i-1}$ to $x_i$ is unique
and hence the lift map is well-defined
if the geodesic distance of $x_{i-1}$ and $x_i$ is less than $\pi$:
we have $c_{d_0 d_1}(t) = \exp_{d_0}(t \log_{d_0}(d_1))$ on $S^2$ and
$c_{R_0 R_1}(t) = \exp_{R_0}(t \log_{R_0}(R_1)) = R_0 \exp(t \log(R_0\tp R_1))$
on \SO3.
The difficulty is that we have to decide numerically
whether the geodesic segment crosses a singularity
(the lifted curve crosses $\bd A_k$ in $T_dS^2$ or $\bd S_k$ in \so3)
in order to determine the proper branch of $\log$ for $y_i$.

All numerical issues discussed below refer to
double precision floating point arithmetic
(64~bit numbers in IEEE 754 format).

\subsection{Directors and Rotations}
\label{sec:num:dir+rot}

Starting with $S^2$, we implement the exponential map $\exp_d$
essentially according to \eqref{eq:s2:expd},
\begin{align}
  \label{eq:s2:expd-num}
  \theta &\define \enorm{v}, &
  \exp_d(v) &=
  \begin{cases}
    d, & \theta \le \epsM, \\
    \displaystyle
    (\cos \theta) d + \frac{\sin \theta}{\theta} v, & \theta > \epsM.
  \end{cases}
\end{align}
Note that we use the machine precision $\epsM > 0$
as threshold in the case distinction
to guard against roundoff errors.
For $\epsM = 0$ (exact arithmetic)
this gives precisely definition \eqref{eq:s2:expd}.
The logarithm map is implemented according to \eqref{eq:s2:logd}
with a second threshold $\delta \in (0, 1)$,
\begin{align}
  \label{eq:s2:logd-num}
  \begin{aligned}
    c &\define \sprod{d}{e}, \\
    p &\define e - c d, \\
    s &\define \enorm{p},
  \end{aligned}
  && \log_d(e) &=
  \begin{cases}
    0, & s \le \epsM, \\
    (\frac1s \arccos c) p, & s > \epsM, \ c > -1 + \delta, \\
    (\arccos c) (\frac1s p), & s > \epsM, \ c \le -1 + \delta.
  \end{cases}
\end{align}
Note that we have $c = \cos \theta$ and $s = \sin \theta$
where $\theta \in (0, \pi]$ is the angle between $d$ and~$e$,
which implies $1/s \to \oo$ for $\theta \to 0$ and for $\theta \to \pi$.
The evaluation order in case 3 would suffice to avoid the singularity
in both cases since $\frac1s p \in S^2$.
However, for $\theta \gtrapprox 0$ we prefer case 2
as it saves two divisions while remaining stable:
$\frac1s \arccos c = \theta / \sin \theta \approx 1$.
Hence the value of $\delta$ is irrelevant in exact arithmetic,
and for $\epsM = 0$ we obtain the original definition \eqref{eq:s2:logd}.
Implementing the lift from $S^2$ into $T_dS^2$ is less straightforward
and will be discussed in \cref{sec:num:lifts}.

How do we determine the threshold value $\delta$?
Our goal is to use case 3 only when necessary
while keeping the \emph{lift-and-project} error
$\norm{\exp_d(\log_d(e)) - e}_\oo$
as small as possible for all $d, e \in S^2$.
To this end, we perform numerical experiments with
\num{9260592} directors in \num{3086864} rotation matrices
from the four \SO3 examples in \cref{sec:examples}.
We pick the smallest negative power of 2 such that
case 3 is consistently more accurate than case 2
for $c \in [-1, -1 + \delta]$.
This yields the value
\begin{equation*}
  \delta = 2^{-21} = 2^{31} \epsM = 32 \sqrt{\epsM}.
\end{equation*}

In \SO3 we implement the exponential map
according to \eqref{eq:so3:exp} and \eqref{eq:so3:exp-aa},
\begin{align}
  \label{eq:so3:exp-num}
  \begin{aligned}
    \theta &\define \norm{N}, \\
    U &\define N / \theta,
  \end{aligned}
  &&
  \exp(N) &=
  \begin{cases}
    I, & \theta \le \epsM, \\
    I + (\sin \theta) U +
    (1 - \cos \theta) U^2, & \theta > \epsM.
  \end{cases}
\end{align}
As in \eqref{eq:s2:expd-num}, $\epsM > 0$ guards against roundoff errors
and $\epsM = 0$ recovers the exact definition.
Of  course, the axis matrix $U$ is only computed in case 2.

The logarithm map \eqref{eq:so3:log}
requires again a more careful case distinction
for $\theta \gtrapprox 0$ and
for $\theta \lessapprox \pi$.
In the case $\theta \gtrapprox 0$ we set $\log(R) = 0$
if $\theta$ is extremely small, otherwise we use the approximation
$R = \exp(N) \approx I + N$, giving
$N = \frac12 (N - N\tp) \approx \frac12 (R - R\tp)$.
In the case $\theta \lessapprox \pi$ we also need to choose the sign of $U$,
\ie, one of the two possible directions of the axis.
We use the approximation
$R \approx \exp(\pi U) = I + 2 U^2 = 2 u u\tp - I$
to determine $\abs{u}$ from the \emph{diagonal} entries of $R$
and the signs of $u_i$ from the \emph{off-diagonal} entries, as follows.
Set $u_i \define (\frac12 R_{ii} + \frac12)^{1/2}$.
If $R_{21} < 0$, flip $u_2$ to $-u_2$.
If $R_{31} < 0$ or $R_{32} < 0 \wedge u_1 = 0$, flip $u_3$ to $-u_3$.
Call the result $u(R)$. Its first nonzero entry will be positive,
that is, the axis direction is fixed
by requiring $u(R)$ to be \emph{lexicographically positive}.
Now implement $\log$ as follows with thresholds
$\delta_1, \delta_2, \delta_3 \in (0, 1)$,
\begin{align}
  \label{eq:so3:log-num}
  \begin{aligned}
    c &\define \frac{\Tr(R) - 1}{2}, \\
    \theta &\define \arccos c,
  \end{aligned}
  &&
  \log(R)
  &=
  \begin{cases}
    0, & c \in [1 - \delta_1, 1], \\
    \frac12 (R - R\tp), & c \in [1 - \delta_2, 1 - \delta_1), \\
    \frac12 (\theta / \sin \theta) (R - R\tp), &
    c \in (-1 + \delta_3, 1 - \delta_2), \\
    \pi \Skew(u(R)), & c \in [-1, -1 + \delta_3].
  \end{cases}
\end{align}
Here $\theta$ is only computed in case 3,
and setting $\delta_1 = \delta_2 = \delta_3 = 0$
recovers definition \eqref{eq:so3:log}.
Again the threshold values are chosen
as the smallest negative powers of $2$ that minimize the
lift-and-project error $\mnorm{\exp(\log(R)) - R}$,
where $\mnorm{}$ on $\R^{3 \x 3}$ is $\norm[\oo]{}$ on $\R^9$.
We obtain the values
\begin{align*}
  \delta_1 &= 2^{-50} = 4 \epsM, &
  \delta_2 &= 2^{-42} = 1024 \epsM, &
  \delta_3 &= 2^{-28} = 2^{24} \eps_M = \frac14 \sqrt{\epsM}.
\end{align*}
In the numerical experiments we use the same
\num{3086864} rotation matrices as above; see \cref{sec:examples}.
Considering the value $c = \frac12 (\Tr(R) - 1)$ of these matrices,
$59$ are in the range $[1 - 16 \epsM, 1)$ relevant for $\delta_1$,
$278$ are in the range $(1 - 4096 \epsM, 1 - 16 \epsM)$ relevant for $\delta_2$,
and $68$ are in the range $(-1, -1 + \sqrt{\epsM}]$ relevant for $\delta_3$.

\subsection{Lifts and Singularities}
\label{sec:num:lifts}

To construct $\lift_{d, \theta u}(e)$ for $d, e \in S^2$
and $\theta u \in \_A_k \subset T_dS^2$ numerically,
we first check whether $e = \pm d$ holds numerically
to return a vector on the appropriate part of $\bd A_k$;
see \cref{alg:s2:liftmap}.
Otherwise we compute the quantities $c, p, s$ and the angle $\arccos c$
of $\log_d(e)$ from \eqref{eq:s2:logd-num}
to determine the proper branch $\ell$ for
$\lift_{d, \theta u}(e) = \log_{d,\ell}(e) = \Phi_{0\ell}(\log_d(e))$
and return the result.
Here the axis-angle decomposition $\theta u$ is unique
except in the case $\theta = 0$ where $u$ is arbitrary
and we simply choose $u = 0$.

The lift of a sequence of elements $x_0, \dots, x_n \in M$
into $T_{x_0}M$ is given in \cref{alg:s2:lift}.
In the algorithms we use the parameters $\tau_1 = 8\epsM$
(the current snapshot is numerically identical to the last one),
$\tau_2 = \epsM$ (geodesic passes through $\pm d$),
and $\tau_3 = 4\epsM$
(geodesic is close to one of the two singularities at $\pm d$).
All parameters have been determined by numerical experiments.

\begin{algorithm}[t]
  \caption{Lift map for $S^2$: $\lift_{d, \theta u}(e)$.}
  \label{alg:s2:liftmap}
  \SetKwInOut{Input}{Input}
  \SetKwInOut{Output}{Output}
  \def\odd{\textup{odd}}
  \Input{\ directors $d, e \in S^2$,
    tangent vector $\theta u \in \_A_k \subset T_dS^2$ in axis-angle form}
  \Output{\ lift $v \in T_dS^2$ of $e$ with respect to $\theta u$}
  \BlankLine
  $c \define \sprod{d}{e}$\;
  $\odd \define k \bmod 2$ (0 or 1)\;
  \If({($e \approx +d$:
    choose $v$ on ``even'' component of $\bd A_k$)}){$1 - c \le \tau_3$}{
    \Return $v \define (k + \odd) \pi u$\;
  }
  \If({($e \approx -d$:
    choose $v$ on ``odd'' component of $\bd A_k$)}){$1 + c \le \tau_3$}{
    \Return $v \define (k + 1 - \odd) \pi u$\;
  }
  $\ell \define k$\;
  $\phi \define \arccos c$\;
  $p \define e - c d$, $w \define (-1)^k p / \enorm{p}$\;
  \If({(geodesic
    passes through $\pm d$)}){
    $\sprod{u}{w} + 1 \le \tau_2$}{
    $w \define -w$\;
    \eIf{$\enorm{((\ell + 2 - \odd) \pi - (-1)^k \phi) w - \theta u} < \pi$
      and ($k > 0$ or $\phi > 1$)}{
      $\texttt{++}\ell$\;
    }{
        \eIf({(prevent change to negative branch)}){$\ell == 0$}{
            $w \define -w$\;
        }{
            $\texttt{--}\ell$\;
        }
    }
    $\odd \define \ell \bmod 2$\;
  }
  $\phi \define (\ell + \odd) \pi + (-1)^k \phi$\;
  \Return $v \define \phi w$\;
\end{algorithm}
\begin{algorithm}[t]
  \caption{Lift of discrete snapshots from a manifold $M$.}
  \label{alg:s2:lift}
  \SetKwInOut{Input}{Input}
  \SetKwInOut{Output}{Output}
  \Input{\ samples $x_0, \dots, x_n \in M$}
  \Output{\ lift $v_0, \dots, v_n \in T_{x_0}M$}
  \BlankLine
  $v_0 \define 0$\;
  \For{$i = 1, \dots, n$}{
    \If({($x_i, x_{i-1}$ are numerically identical)}){
      $\operatorname{dist}(x_i, x_{i-1}) \le \tau_1$}{
      $v_i \define v_{i-1}$\;
    }
    \Else{$v_i \define \lift_{x_0,v_{i-1}} (x_i)$\;}
  }
\end{algorithm}

To construct the lift map for \SO3 numerically,
we compute the numerical logarithm map \eqref{eq:so3:log-num}
in the following axis-angle form with
$c = \frac12 (\Tr(R) - 1)$, $V \define R - R\tp$:
\begin{align*}
  \theta &= 0,            & U &= 0,                      &\text{case 1}, \\
  \theta &= \norm{V} / 2, & U &= V / \theta,             &\text{case 2}, \\
  \theta &= \arccos c,    & U &= V / (2 \sqrt{1 - c^2}), &\text{case 3}, \\
  \theta &= \pi,          & U &= \Skew(u(R)),            &\text{case 4}.
\end{align*}
Again $U$ is arbitrary in the case $\theta = 0$ and we have chosen $U = 0$.
Then, given $N \in \so3$, we need to determine the proper branch $\ell$
for $\lift_N(R) = \log_\ell(R) = \Phi_{0\ell}(\theta U)$
where we distinguish three cases, see \cref{alg:so3:liftmap}:
\emph{i)} if $\norm{N} + \theta < \pi$, then $\ell = 0$ (quick return);
\emph{ii)} if $N$ and $U$ are linearly dependent,
then $\exp(N)$ and $R$ are rotations about a common axis
and the computation of $\ell$ simplifies; and
\emph{iii)} otherwise we have the slightly more complicated general case.
Note that ``even'' branch boundaries at angles $2 k \pi$
(red in \cref{fig:shells})
can only be crossed in case \emph{ii)}.
In \cref{alg:so3:liftmap} we use the parameter $\epsO3 = 8 \epsM$,
which is the best accuracy that can be guaranteed for the
numerical representation of orthogonal matrices
with respect to the \emph{orthogonality error}
$\orterr(R) \define \mnorm{R\tp R - I}$,
as determined by numerical experiments; \cf \cref{sec:accuracy}.

\begin{algorithm}[t]
  \caption{Lift map for \SO3: $\lift_N(R)$.}
  \label{alg:so3:liftmap}
  \SetKwInOut{Input}{Input}
  \SetKwInOut{Output}{Output}
  \def\odd{\textup{odd}}
  \Input{\ tangent matrix $N \in \so3$, rotation matrix $R \in \SO3$}
  \Output{\ lift $\theta U$ of $R$ with respect to $N$}
  \BlankLine
  $\theta_0 \define \norm N$\;
  $\theta U \define \log(R)$
  (initial axis-angle-form: a bijection $\Phi_{0\ell}$ may be applied later)\;
  \If{$\theta_0 + \theta < \pi$}{
    \Return $\theta U$
  }
  $s \define \sprod{N}{U}$\;
  \If({($N, U$ are numerically linearly dependent)})
    {$\Abs{\abs{s} - \theta_0\strut} < 8 \epsO3 \theta_0$}{
    \lIf{$s < 0$}{
      $\_\theta \define -\theta$%
    }
    \lElse{
      $\_\theta \define +\theta$%
    }
    \If{$\theta_0 - \_\theta < \pi$}{
      \Return $\theta U$
    }
    $\ell \define \floor{(\theta_0 - \_\theta) / \pi}$\;
    \lIf{$\ell$ is odd}{$\texttt{++}\ell$}
    \If{$\_\theta > 0$}{
      \Return $+(\_\theta + \ell \pi) U$
    }
    \Else{
      \Return $-(\_\theta + \ell \pi) U$
    }
  }
  $k \define \ell \define \floor{\theta_0 / \pi}$\;
  $\odd \define k \bmod 2$\;
  \If{$k \ne 0$}{$N \define [1 - (k + \odd) \pi / \theta_0] N$}
  \If{$\norm{\theta U - N} > \norm{(\theta - 2 \pi) U - N}$}{
    \lIf{\odd}{
      $\texttt{--}\ell$%
    }
    \lElse{
      $\texttt{++}\ell$%
    }
  }
  \lIf{$\ell$ is odd}{
    $\theta \define \theta - (\ell + 1) \pi$%
  }
  \lElse{
    $\theta \define \theta + \ell \pi$%
  }
  \Return $\theta U$
\end{algorithm}

\subsection{PCA in the Tangent Space}
\label{sec:pca}

In dynamics of hybrid mechanical systems
we typically perform PGA on a \emph{product manifold}
of the form $M = (\R^3 \x \SO3)^{K_1} \x (\R^3 \x S^2)^{K_2} \x (\R^3)^{K_3}$
that results from a combination of rigid bodies
and finite element discretizations of beams, shells,
and 3-dimensional elastic bodies,
or possibly on a submanifold of $M$ arising due to holonomic constraints.
The principal component analysis in the tangent space $T_{x_0}M$
is performed on the \emph{snapshot matrix} obtained from lifting
$x_1, \dots, x_n \in M$ where $x_i = c(t_i)$:
\begin{equation}
  \label{eq:lift:Y}
  Y \define \mat[3]{y_1 & \dots & y_n} \in (T_{x_0}M)^n \cong \R^{m \x n}.
\end{equation}
Then, as usual, we compute the singular value decomposition (SVD)
\begin{align*}
  Y &= U \Sigma V\tp = \sum_{j = 1}^r \sigma_j u_j v_j\tp,
  & r &= \rank Y \le \min(m, n),
\end{align*}
where $U \in O(m)$ and $V \in O(n)$ are orthogonal matrices
that contain the left and right \emph{singular vectors}
$u_j$ and $v_j$, respectively,
and $\Sigma \in \R^{m \x n}$ is a rectangular diagonal matrix
that contains the \emph{singular values}
$\sigma_1 \ge \dotsm \ge \sigma_r > 0$.
Best rank $p$ approximations of $Y$ with respect to both
the spectral norm $\enorm{}$ and the Frobenius norm $\norm{}_F$
are then obtained for each $p \le r$ by the \emph{truncated SVD}
\begin{equation}
  \label{eq:trsvd}
  Y_p = U \Sigma_p V\tp = \sum_{j = 1}^p \sigma_j u_j v_j\tp.
\end{equation}
The matrices $\sigma_j u_j v_j\tp$ in this sum
are the $p$ principal components (or principal modes)
of the lifted discrete trajectory represented by $Y$.
The resulting approximated snapshots in $M$ are finally obtained
by mapping the columns $y_{pk}$ of the matrix $Y_p$ back to $M$,
\begin{align}
  \label{eq:tr-proj}
  x_{pk} &= c_p(t_k) = \exp_{x_0}(y_{pk}), &
  y_{pk} &= Y_p e_k = \sum_{j = 1}^p \sigma_j u_j v_{jk}, &
  k &= 1, \dots, n.
\end{align}
Note that the reference point $x_0$ does not change here:
it is the given initial point of the discrete trajectory,
and it always maps to $y_0 = 0 \in T_{x_0}M$
which we have therefore dropped in forming~$Y$.

The situation is different when considering multiple trajectories
of a dynamic system simultaneously.
In that case we compute a (not necessarily unique) \emph{intrinsic mean} $\_x$
of all points on all $k$ discrete trajectories \cite{Fletcher2004},
\begin{equation*}
  \sum_{i=1}^k \sum_{j = 0}^{n_i} d_M(x_{ij}, \_x)^2 \stackrel!=
  \min_{x \in M} \sum_{i=1}^k \sum_{j = 0}^{n_i} d_M(x_{ij}, x)^2,
\end{equation*}
we do not drop any initial points, and the snapshot matrix contains
the entire set of concatenated trajectories
lifted into the tangent space at the mean $\_x$:
\begin{align*}
  Y &\define \mat[7]{y_{10} & \dots & y_{1n_1} & \dots & y_{k0} & \dots & y_{kn_k}}
  \in (T_{\_x}M)^n \cong \R^{m \x n}, \\
  n &\define (n_1 + 1) + \dotsm + (n_k + 1).
\end{align*}
The intrinsic mean is typically computed by a gradient descent algorithm
with fixed step length $\alpha$
where the iterate $\_x_\ell$ is updated by the rule
\begin{equation*}
  \_x_{\ell + 1} = \exp_{\_x_\ell} \left(
    \frac{\alpha}{n} \sum\nolimits_{i, j} \log_{\_x_\ell}(x_{ij}) \right).
\end{equation*}
Of course, on \SO3 we always lift into \so3
and obtain the required tangent vectors
by simple multiplication with the intrinsic mean $\_R$
or with the initial rotation $R_0$, if desired.
Here the update rule for the intrinsic mean of matrices $R_{ij}$ becomes
\begin{equation*}
  \_R_{\ell + 1} = \_R_\ell \exp \left(
    \frac{\alpha}{n} \sum\nolimits_{i, j} \log(\_R_\ell\tp R_{ij}) \right).
\end{equation*}

\subsection{Discussion of Accuracy}
\label{sec:accuracy}

In this section we discuss the numerical accuracy of
the nonlinear mappings implemented in \cref{sec:num}.
We start with the projections
\eqref{eq:so3:proj} for \SO3 and \eqref{eq:s2:proj} for $S^2$,
\ie, the accuracy of representing elements of these manifolds.
Then we address the lift-and-project mappings
\explog and \explift for \SO3 and $S^2$.
For \SO3, error statistics are based on \num{3086864} rotation matrices
from the four examples in \cref{sec:examples}.
For $S^2$, we have \num{655968} directors
from the three $S^2$ examples in \cref{sec:examples}
plus \num{9260592} directors from the four \SO3 examples.

As already mentioned in \cref{sec:num:lifts},
we can bound the orthogonality error on \SO3,
$\orterr(R) = \mnorm{R\tp R - I}$, by \epsO3.
If a rotation matrix $R$
(from a numerical simulation, for instance)
does not satisfy that bound,
we replace it with its projection onto \SO3
according to \eqref{eq:so3:proj}
using the \LAPACK routine \texttt{dgesvd}
\cite{Anderson_et_al:1992}.
Our numerical experiments show that this gives
an accuracy of $\frac74 \epsO3$ or better.
Then, if necessary, a single Newton iteration for the equation
$R\tp R - I = 0$ reduces the error to at most \epsO3.

For $S^2$ the situation is even better:
when applying the projection \eqref{eq:s2:proj},
the \emph{normalization error} $\abs{1 - \enorm{d}}$ is zero in most cases
(oscillating beam 79.2\%, triple pendulum 77.2\%,
rubber rod 80.0\%, wind turbine 64.2\%,
shell pendulum 63.1\%, free flying plate 61.0\%, tumbling cylinder 62.4\%),
and it lies in the range $[\frac12 \epsM, 2 \epsM]$ in all other cases.

Orthogonality errors $\orterr(R)$,
$\orterr(\exp(\log(R)))$ and $\orterr(\exp(\lift(R)))$
on \SO3 are presented in \cref{fig:so3:bin-error}
for the triple pendulum and the wind turbine.
(The oscillating beam and the rubber rod
perform 2d motions like the triple pendulum;
they behave similarly and they are described in the text,
but we omit visualizations.)
Input matrices ($R$) are projected onto \SO3 as just described
if the initial error exceeds $\epsO3 = 8 \epsM$.
This applies to \num{600111} matrices for the wind turbine
and none for the other cases.
The maximal output error is $15 \epsM$ for the wind turbine
and $\epsO3$ or better in all other cases.
The number of errors exceeding \epsO3 is
2528 or $0.39\%$ for \explog and
2492 or $0.38\%$ for \explift.
As explained, a single Newton iteration
will reduce all errors to at most \epsO3 if desired.
In summary, our results show that \explog
and \explift map \SO3 indeed to itself
up to unavoidable numerical roundoff errors.

\pgfplotsset{/pgfplots/tick align=outside,/pgfplots/tick pos=left,
  /pgfplots/bar cycle list/.style={/pgfplots/cycle list={
      {AOblue,fill=AOblue!30,mark=none},
      {AOred,fill=AOred!30,mark=none},
      {AOgreen!75!black,fill=AOgreen!36,mark=none},
    }},/pgfplots/ybar=0.4pt}
\def\PTitle{\node at (rel axis cs:0.5,0.87)}

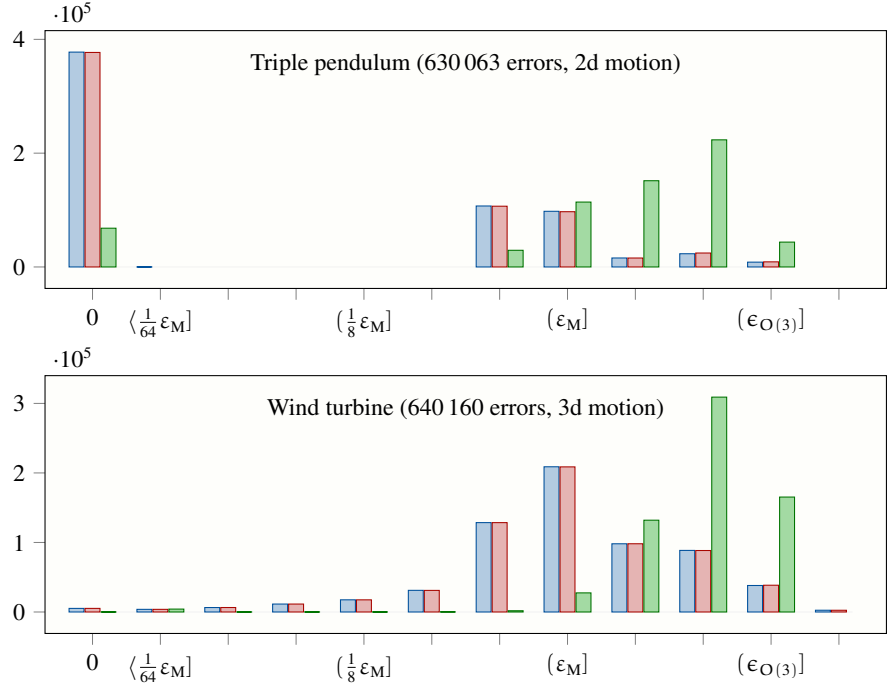
\begin{figure}[tp]
  \small
  \begin{tikzpicture}
    \begin{axis}[width=\textwidth,height=5cm,xmin=-0.7,xmax=11.7,
      bar width=0.22,xtick distance=1,xticklabels={,0,
        $\langle\frac1{64}\epsM]$,,,$(\frac18\epsM]$,,,$(\epsM]$,,,$(\epsO3]$}
      ]
      \draw[black!5] (-0.2,0) -- (10.2,0);
      \addplot coordinates {
        (0,377584) (1,5) (6,107157) (7,97830) (8,15785) (9,23237) (10,8465)
      };
      \addplot coordinates {
        (0,376987) (6,106763) (7,97109) (8,15720) (9,24504) (10,8980)
      };
      \addplot coordinates {
        (0,68109) (6,29298) (7,114036) (8,151557) (9,223394) (10,43669)
      };
      \PTitle{Triple pendulum (\num{630063} errors, 2d motion)};
    \end{axis}
  \end{tikzpicture}
  \begin{tikzpicture}
    \begin{axis}[width=\textwidth,height=5cm,xmin=-0.7,xmax=11.7,
      bar width=0.22,xtick distance=1,xticklabels={,0,
        $\langle\frac1{64}\epsM]$,,,$(\frac18\epsM]$,,,$(\epsM]$,,,$(\epsO3]$}
      ]
      \draw[black!5] (-0.2,0) -- (11.2,0);
      \addplot coordinates {
        (0,5169) (1,3806) (2,6415) (3,11420) (4,17499) (5,31134) (6,128547) (7,208856) (8,98108) (9,88569) (10,38109) (11,2528)
      };
      \addplot coordinates {
        (0,5183) (1,3809) (2,6439) (3,11416) (4,17470) (5,31058) (6,128531) (7,208707) (8,98135) (9,88422) (10,38498) (11,2492)
      };
      \addplot coordinates {
        (0,22) (1,4167) (2,33) (3,53) (4,72) (5,188) (6,1792) (7,27479) (8,132027) (9,309026) (10,165301)
      };
      \PTitle{Wind turbine (\num{640160} errors, 3d motion)};
    \end{axis}
  \end{tikzpicture}
  \caption{Bin counts of orthogonality error:
    $R$ (green), $\exp(\log(R))$ (blue), $\exp(\lift(R))$ (red).
    Error bins: $\langle x] = (0, x]$, $(x] = (\frac12 x, x]$.}
  \label{fig:so3:bin-error}
\end{figure}

No visual or tabular presentations are given for the
normality errors on $S^2$ since all of them are negligible.

\begin{figure}[tp]
  \small
  \begin{tikzpicture}
    \begin{axis}[width=\textwidth,height=5cm,xmin=-1.5,xmax=42.5,
      bar width=0.25,xtick={0,...,41},
      xticklabels={0,*,\rlap{\llap{$($}$0]$},,,,$(4]$,,,,,$(9]$,,,,,
        $(14]$,,,,,$(19]$,,,,,$(24]$,,,,,$(29]$,,,,,$(34]$,,,,,$(39]$}
      ]
      \draw[black!5] (-0.2,0) -- (41.2,0);
      \addplot coordinates {
        (0,11653) (1,84999) (2,81849) (3,142009) (4,99913) (5,46350) (6,37386) (7,35424) (8,26212) (9,18204) (10,13013) (11,8881) (12,6694) (13,4845) (14,3933) (15,2643) (16,1739) (17,1254) (18,853) (19,563) (20,425) (21,306) (22,230) (23,184) (24,117) (25,88) (26,72) (27,60) (28,65) (29,48) (30,17) (31,8) (32,1) (38,2) (39,6) (40,7) (41,10)
      };
      \addplot coordinates {
        (0,11581) (1,84691) (2,81849) (3,141761) (4,99692) (5,47121) (6,38414) (7,36384) (8,28105) (9,22010) (10,20520) (11,15566) (12,1682) (13,326) (14,29) (15,18) (16,13) (17,16) (18,17) (19,15) (20,25) (21,24) (22,26) (23,24) (24,19) (25,23) (26,18) (27,23) (28,35) (29,31) (30,5)
      };
      \PTitle{Triple pendulum (\num{630063} errors, 2d motion)};
    \end{axis}
  \end{tikzpicture}
  \begin{tikzpicture}
    \begin{axis}[width=\textwidth,height=5cm,xmin=-1.5,xmax=42.5,
      bar width=0.25,xtick={0,...,41},
      xticklabels={0,*,\rlap{\llap{$($}$0]$},,,,$(4]$,,,,,$(9]$,,,,,
        $(14]$,,,,,$(19]$,,,,,$(24]$,,,,,$(29]$,,,,,$(34]$,,,,,$(39]$}
      ]
      \draw[black!5] (-0.2,0) -- (41.2,0);
      \addplot[AOblue,fill=AOblue!30] coordinates {
        (0,14) (1,12138) (2,51699) (3,176109) (4,201322) (5,87632) (6,31884) (7,24997) (8,19345) (9,13084) (10,8671) (11,5169) (12,2859) (13,1466) (14,1009) (15,832) (16,554) (17,383) (18,265) (19,166) (20,124) (21,75) (22,69) (23,42) (24,36) (25,35) (26,35) (27,41) (28,45) (29,43) (30,6) (31,1) (35,1) (40,4) (41,5)
      };
      \addplot[AOred,fill=AOred!30] coordinates {
        (0,14) (1,11926) (2,49891) (3,167056) (4,189708) (5,95395) (6,45082) (7,26762) (8,19325) (9,13069) (10,8665) (11,5171) (12,2857) (13,1469) (14,1008) (15,832) (16,554) (17,383) (18,265) (19,166) (20,124) (21,75) (22,69) (23,42) (24,36) (25,35) (26,35) (27,41) (28,45) (29,43) (30,6) (31,1) (35,1) (40,4) (41,5)
      };
      \PTitle{Wind turbine (\num{640160} errors, 3d motion)};
    \end{axis}
  \end{tikzpicture}
  \caption{Bin counts of \SO3 lift-and-project error:
    \explog (blue), \explift (red).
    Error bins: $* = (0, \frac12 \epsM]$, $(n] = (2^{n-1} \epsM, 2^n \epsM]$.}
  \label{fig:so3:bin-dist}
\end{figure}
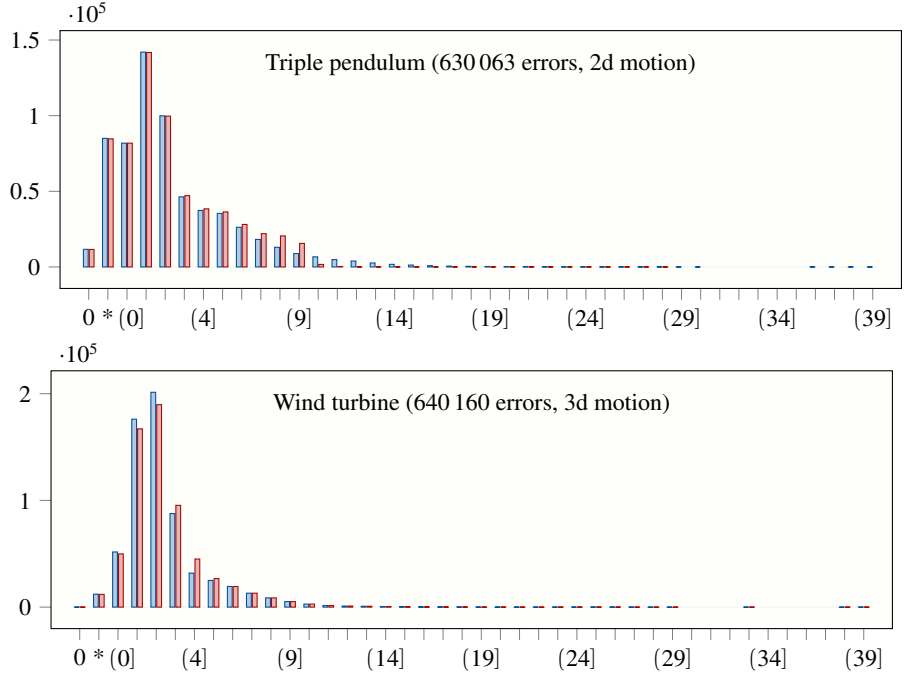

The lift-and-project errors
$\mnorm{\exp(\log(R)) - R}$ and
$\mnorm{\exp(\lift(R)) - R}$
for \SO3 are shown in \cref{fig:so3:bin-dist}.
These are crucial for PGA.
We observe that the vast majority of errors
are below $2^{12} \epsM \approx \num{2e-12}$,
but a tiny fraction range up to almost \num{1e-4}.
More detailed information is provided in \cref{tab:so3:dist}.
A glance at the angular distribution of the
lift-and-project errors in \cref{fig:so3:ang-dist}
reveals that, as expected, all these large errors
occur at rotation angles close to zero and $\pi$:
they are caused by the singularities of the logarithm map.

\begin{table}[tp]
  \small
  \caption{Error statistics of \explog and \explift on \SO3.}
  \label{tab:so3:dist}
  \begin{tabular}{llrrlrr}
    \toprule
    Dynamic system &
    \multicolumn3c{error of \explog} & \multicolumn3c{error of \explift} \\
    & maximal & \multicolumn2r{large ($> 2^{12} \epsM$)}
    & maximal & \multicolumn2r{large ($> 2^{12} \epsM$)} \\
    \midrule
    Oscillating beam &
    \num{3.9e-13} &    0 &    0\% & \num{2.3e-13} &    0 &    0\% \\
    Triple pendulum &
    \num{8.5e-05} & 8698 & 1.38\% & \num{3.1e-08} &  332 & 0.05\% \\
    Rubber rod &
    \num{9.5e-08} &   59 & 1.13\% & \num{1.5e-08} &   11 & 0.21\% \\
    Wind turbine &
    \num{8.4e-05} & 2762 & 0.43\% & \num{8.4e-05} & 2762 & 0.43\% \\
    \bottomrule
  \end{tabular}
\end{table}

\newcommand\errimg[4][35]{\twtrimg[0.49]
  {18 20 46 #1}{#2.png}%
  \llap{\vtop to 0pt{\vss
      \hbox to0.49\textwidth{\qquad\hss#3\hss}%
      \hbox to0.49\textwidth{\qquad\hss(\num{#4} errors)\hss
      }\kern0.255\textwidth}}}%
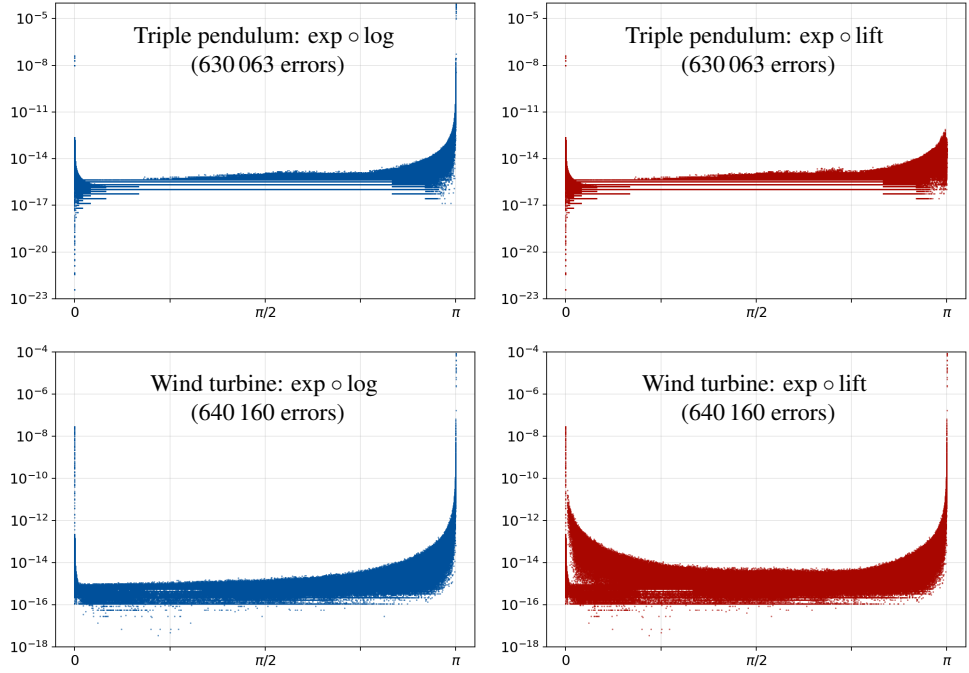
\begin{figure}[tp]
  \small
  \ifcase0
  \errimg[41]{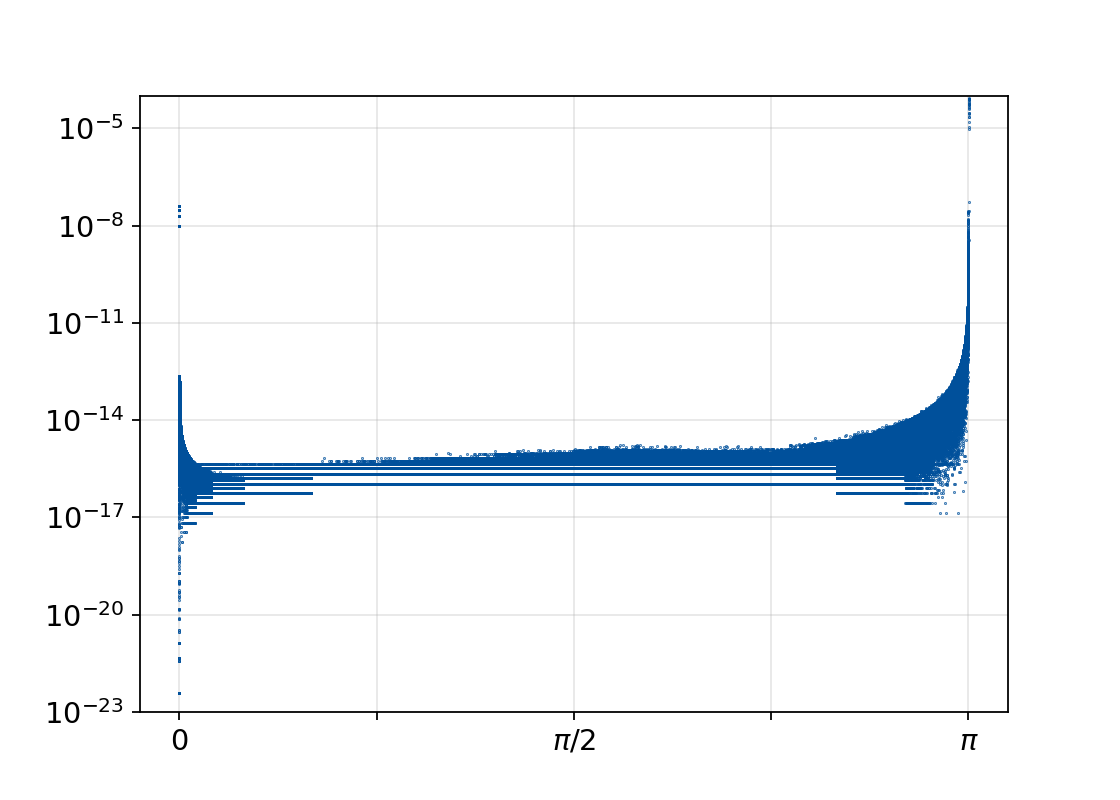}{Triple pendulum: \explog}{630063}\hfill
  \errimg[41]{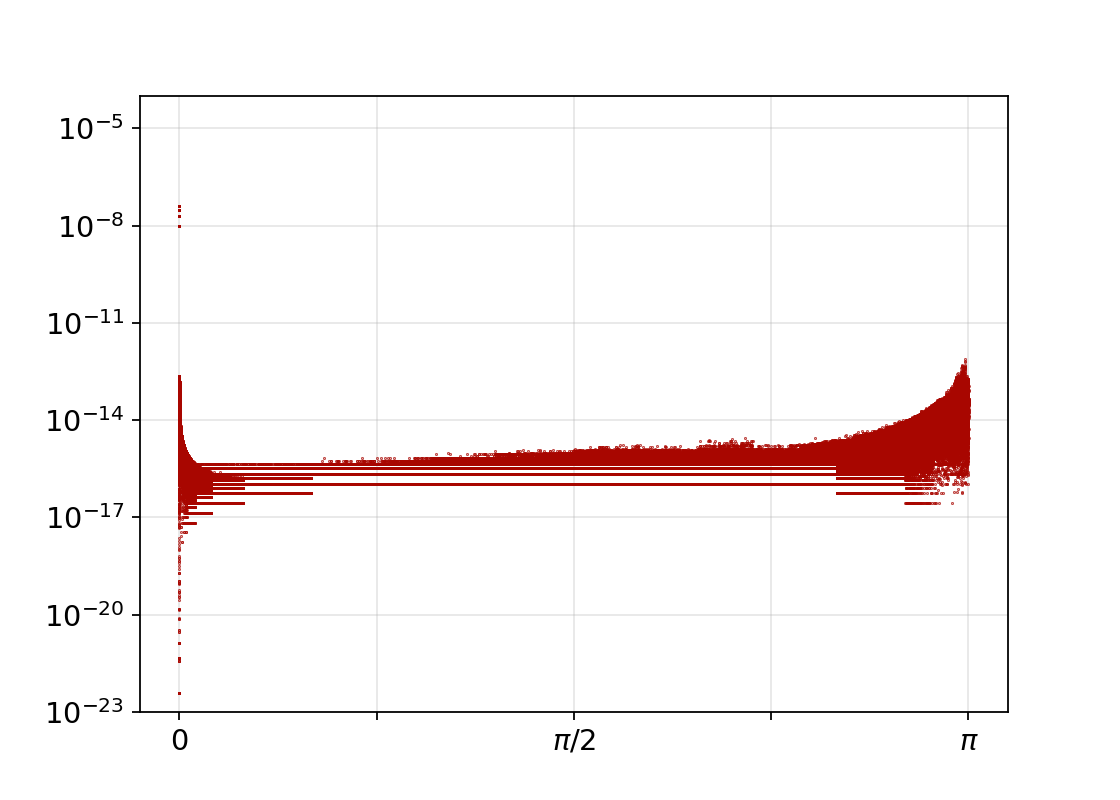}{Triple pendulum: \explift}{630063}\\[2ex]
  \errimg{wind_log}{Wind turbine: \explog}{640160}\hfill
  \errimg{wind_lift}{Wind turbine: \explift}{640160}%
  \or
  \begin{tikzpicture}
    \begin{semilogyaxis}[width=0.49\textwidth,height=5cm,xmin=-0.2,xmax=3.3416,
        xtick distance=pi/4,xticklabels={,0,,,,$\pi$}
      ]
      \addplot[AOblue,only marks,mark=*,mark size=0.4pt]
      table[x index=0,y index=2] {tmp/pend_log};
      \PTitle{Triple pendulum};
    \end{semilogyaxis}
  \end{tikzpicture}
  \begin{tikzpicture}
    \begin{semilogyaxis}[width=0.49\textwidth,height=5cm,xmin=-0.2,xmax=3.3416,
        xtick distance=pi/4,xticklabels={,0,,,,$\pi$}
      ]
      \addplot[AOblue,only marks,mark=*,mark size=0.4pt]
      table[x index=0,y index=2] {tmp/wind_log};
      \PTitle{Wind turbine};
    \end{semilogyaxis}
  \end{tikzpicture}
  \begin{tikzpicture}
    \begin{semilogyaxis}[width=0.49\textwidth,height=5cm,xmin=-0.2,xmax=3.3416,
        xtick distance=pi/4,xticklabels={,0,,,,$\pi$}
      ]
      \addplot[AOred,only marks,mark=*,mark size=0.4pt]
      table[x index=2,y index=4] {tmp/pend_lift};
      \PTitle{Triple pendulum};
    \end{semilogyaxis}
  \end{tikzpicture}
  \hfil
  \begin{tikzpicture}
    \begin{semilogyaxis}[width=0.49\textwidth,height=5cm,xmin=-0.2,xmax=3.3416,
        xtick distance=pi/4,xticklabels={,0,,,,$\pi$}
      ]
      \addplot[AOred,only marks,mark=*,mark size=0.4pt]
      table[x index=2,y index=4] {tmp/wind_lift};
      \PTitle{Wind turbine};
    \end{semilogyaxis}
  \end{tikzpicture}
  \fi
  \caption{Angular distributions of \SO3 lift-and-project error.}
  \label{fig:so3:ang-dist}
\end{figure}

The lift-and-project errors for $S^2$ are illustrated
in \cref{fig:s2:bin-dist1} for the three shell examples,
and in \cref{fig:s2:bin-dist3} for the individual directors
of the triple pendulum and the wind turbine.
Again the majority of errors are below $2^{12} \epsM \approx \num{2e-12}$
while the remaining ones are smaller than for \SO3,
with maximum at \num{5e-8}.

\begin{figure}[tp]
  \centering
  \small
  \begin{tikzpicture}
    \begin{axis}[width=\textwidth,height=5cm,xmin=-0.9,xmax=19.9,
      bar width=0.35,xtick distance=1,
      xticklabels={,0,$\langle0]$,,,$(3]$,,,,,$(8]$,,,,,$(13]$,,,,,$(18]$}
      ]
      \draw[black!5] (-0.2,0) -- (19.2,0);
      \addplot coordinates {
        (0,802) (1,12817) (2,3304) (3,1093) (4,618) (5,467) (6,317) (7,297) (8,236) (9,164) (10,130) (11,87) (12,76) (13,44) (14,35) (15,20) (16,17) (17,8) (18,3) (19,1)
      };
      \addplot coordinates {
        (0,597) (1,12430) (2,3819) (3,1174) (4,614) (5,467) (6,317) (7,297) (8,236) (9,164) (10,130) (11,87) (12,76) (13,44) (14,35) (15,20) (16,17) (17,8) (18,3) (19,1)
      };
      \PTitle{Tumbling cylinder (20\,468 errors, 3d motion)};
    \end{axis}
  \end{tikzpicture}
  \begin{tikzpicture}
    \begin{axis}[width=\textwidth,height=5cm,xmin=-0.9,xmax=19.9,
      bar width=0.35,xtick distance=1,
      xticklabels={,0,$\langle0]$,,,$(3]$,,,,,$(8]$,,,,,$(13]$,,,,,$(18]$}
      ]
      \draw[black!5] (-0.2,0) -- (19.2,0);
      \addplot coordinates {
        (0,7928) (1,110930) (2,26637) (3,6987) (4,1305) (5,481) (6,218) (7,141) (8,113) (9,94) (10,75) (11,53) (12,58) (13,41) (14,31) (15,18) (16,20) (17,13) (18,11) (19,1)
      };
      \addplot coordinates {
        (0,6642) (1,107778) (2,30356) (3,7637) (4,1362) (5,489) (6,222) (7,141) (8,113) (9,94) (10,75) (11,53) (12,58) (13,41) (14,31) (15,18) (16,20) (17,13) (18,11) (19,1)
      };
      \PTitle{Free-flying plate (155\,155 errors, 3d motion)};
    \end{axis}
  \end{tikzpicture}
  \begin{tikzpicture}
    \begin{axis}[width=\textwidth,height=5cm,xmin=-0.9,xmax=29.9,
      bar width=0.35,xtick distance=1,xticklabels={,0,$\langle0]$,,,
        $(3]$,,,,,$(8]$,,,,,$(13]$,,,,,$(18]$,,,,,$(23]$,,,,,$(28]$}
      ]
      \draw[black!5] (-0.2,0) -- (29.2,0);
      \addplot coordinates {
        (0,22796) (1,287032) (2,74110) (3,42837) (4,22664) (5,11949) (6,6620) (7,3842) (8,2479) (9,1645) (10,1199) (11,809) (12,537) (13,347) (14,214) (15,152) (16,202) (17,122) (18,124) (19,97) (20,65) (21,97) (22,37) (23,97) (24,78) (25,79) (26,65) (27,40) (28,10)
      };
      \addplot coordinates {
        (0,20687) (1,276053) (2,84705) (3,45146) (4,22899) (5,11904) (6,6612) (7,3844) (8,2485) (9,1636) (10,1202) (11,808) (12,538) (13,347) (14,214) (15,152) (16,202) (17,122) (18,124) (19,97) (20,65) (21,97) (22,32) (23,92) (24,68) (25,64) (26,50) (27,40) (28,25) (29,35)
      };
      \PTitle{Shell pendulum (480\,345 errors, 2d motion)};
    \end{axis}
  \end{tikzpicture}
  \caption{Bin counts of $S^2$ lift-and-project error:
    \explog[d] (blue), \explift[d] (red).
    Error bins: $\langle n] = (0, 2^n \epsM]$,
    $(n] = (2^{n-1} \epsM, 2^n \epsM]$.}
  \label{fig:s2:bin-dist1}
\end{figure}

\begin{figure}[tp]
  \centering
  \small
  \begin{tikzpicture}
    \begin{axis}[width=\textwidth,height=5cm,xmin=-0.9,xmax=29.9,
      bar width=0.35,xtick distance=1,xticklabels={,0,$\langle0]$,,,
        $(3]$,,,,,$(8]$,,,,,$(13]$,,,,,$(18]$,,,,,$(23]$,,,,,$(28]$}
      ]
      \draw[black!5] (-0.2,0) -- (30.2,0);
      \addplot coordinates {
        (0,719388) (1,715264) (2,189468) (3,120437) (4,61150) (5,33518) (6,19154) (7,10142) (8,6134) (9,4170) (10,2852) (11,1996) (12,1446) (13,1092) (14,812) (15,714) (16,564) (17,410) (18,310) (19,226) (20,184) (21,146) (22,116) (23,112) (24,104) (25,104) (26,86) (27,74) (28,16)
      };
      \addplot coordinates {
        (0,721205) (1,679470) (2,218557) (3,125743) (4,60900) (5,33380) (6,19150) (7,10124) (8,6126) (9,4168) (10,2854) (11,1996) (12,1446) (13,1092) (14,812) (15,714) (16,564) (17,410) (18,310) (19,226) (20,182) (21,144) (22,108) (23,104) (24,84) (25,80) (26,54) (27,70) (28,62) (29,54)
      };
      \PTitle{Triple pendulum with $S^2$ (1\,890\,189 errors, 2d motion)};
    \end{axis}
  \end{tikzpicture}
  \begin{tikzpicture}
    \begin{axis}[width=\textwidth,height=5cm,xmin=-0.9,xmax=29.9,
      bar width=0.35,xtick distance=1,xticklabels={,0,$\langle0]$,,,
        $(3]$,,,,,$(8]$,,,,,$(13]$,,,,,$(18]$,,,,,$(23]$,,,,,$(28]$}
      ]
      \draw[black!5] (-0.2,0) -- (30.2,0);
      \addplot coordinates {
        (0,36093) (1,621166) (2,263162) (3,201579) (4,145890) (5,110108) (6,89661) (7,78944) (8,72797) (9,70012) (10,62062) (11,48571) (12,36023) (13,24789) (14,16579) (15,11764) (16,7736) (17,5320) (18,3990) (19,3241) (20,2525) (21,2052) (22,1588) (23,1295) (24,1025) (25,778) (26,568) (27,748) (28,411) (29,3)
      };
      \addplot coordinates {
        (0,32198) (1,606434) (2,276521) (3,206009) (4,146233) (5,110148) (6,89669) (7,78976) (8,72832) (9,70041) (10,62090) (11,48594) (12,36068) (13,24816) (14,16604) (15,11809) (16,7787) (17,5369) (18,4032) (19,3282) (20,2577) (21,2107) (22,1625) (23,1328) (24,1002) (25,678) (26,428) (27,429) (28,518) (29,276)
      };
      \PTitle{Wind turbine with $S^2$ (1\,920\,480 errors, 2d motion)};
    \end{axis}
  \end{tikzpicture}
  \caption{Bin counts of $S^2$ lift-and-project error for \SO3 directors:
    \explog[d] (blue), \explift[d] (red).
    Error bins: $\langle n] = (0, 2^n \epsM]$,
    $(n] = (2^{n-1} \epsM, 2^n \epsM]$.}
  \label{fig:s2:bin-dist3}
\end{figure}

Excessive lift-and-project errors on \SO3 can of course be reduced
to $2^{12} \epsM$ or even to \epsO3
by performing certain critical operations with higher accuracy
(128~bit quadruple precision might suffice).
This should not even be overly expensive
as it is only required close to the singularities.
Moreover, arbitrary precision software
such as the Multiprecision library of \boost~\cite{Boost}
might be used to save further cost
by avoiding unnecessarily accurate computations.
However, these issues are beyond the scope of this paper,
and we need to point out that the current logarithm and lift maps for \SO3
do not provide sufficient numerical precision
for employing PGA-based model order reduction in dynamic simulations
when matrices in the vicinity of the singularities occur.

For $S^2$ the situation is definitely better,
but with just below eight decimal digits in the worst case
(error \new{\num{4.9e-8}}) it is probably still not good enough
for PGA-based model order reduction in most situations.
Thus, increasing the numerical precision
of the logarithm and lift maps for $S^2$
will be a topic of further research as well.

\begin{figure}
  \small
  \centering
  \errimg{tc_log}{Tumbling cylinder: \explog}{20468}\hfill
  \errimg{tc_lift}{Tumbling cylinder: \explift}{20468}\\[2ex]
  \errimg{ffp_log}{Free-flying plate: \explog}{155155}\hfill
  \errimg{ffp_lift}{Free-flying plate: \explift}{155155}\\[2ex]
  \errimg{sp_log}{Shell pendulum: \explog}{480345}\hfill
  \errimg{sp_lift}{Shell pendulum: \explift}{480345}%
  \caption{Angular distributions of $S^2$ lift-and-project error.}
\end{figure}

\begin{table}[tp]
  \small
  \caption{Error statistics of \explog[d] and \explift[d] on $S^2$.}
  \label{tab:s2:dist}
  \begin{tabular}{llrrlrr}
    \toprule
    Dynamic system &
    \multicolumn3c{error of \explog} & \multicolumn3c{error of \explift} \\
    & maximal & \multicolumn2r{large ($> 2^{12} \epsM$)}
    & maximal & \multicolumn2r{large ($> 2^{12} \epsM$)} \\
    \midrule
    Tumbling cylinder &
    \num{4.0e-11} &   128 & 0.63\% & \num{4.0e-11} &   128 & 0.63\% \\
    Free-flying plate &
    \num{3.6e-11} &   193 & 0.12\% & \num{3.6e-11} &   193 & 0.12\% \\
    Shell pendulum &
    \num{1.7e-08} &  2363 & 0.49\% & \num{4.2e-08} &  2364 & 0.49\% \\
    \midrule
    Triple pendulum &
    \num{1.7e-08} &  5070 & 0.27\% & \num{4.6e-08} &  5070 & 0.27\% \\
    Wind turbine &
    \num{3.3e-08} & 84412 & 4.39\% & \num{4.9e-08} & 84667 & 4.41\% \\
    \bottomrule
  \end{tabular}
\end{table}

\section{Computational Examples}
\label{sec:examples}

\subsection{Examples on \SO3}
\label{sec:examples-SO3}

To carry out our investigations on \SO3, we select four examples:
a swinging rubber rod,
a free-oscillating cantilever beam,
a flexible triple pendulum, and
a horizontal-axis wind turbine.
All four examples have been used for error statistics in \cref{sec:accuracy},
and we have earlier studied them in \cite{Gebhardt_et_al:2019}.
However, for the sake of brevity, the principal geodesic analysis on \SO3
is carried out only for the triple pendulum and the wind turbine.
In addition and due to the director-based formulation adopted in this work,
the principal geodesic analysis on $S^2$ is also carried out
for individual directors of the triple pendulum and wind turbine examples,
whose results are presented in \cref{sec:examples-S2}.

\subsubsection{Rubber rod}

This example considers a swinging rubber rod
that is simply supported at one of its ends
and subject only to the action of its own weight induced by gravity.
Upon sudden release, the rod moves in a vertical plane
exhibiting nonlinear kinematics.
The geometrical and cross-sectional material properties of the rod
are the following:
length \qty{1.0}{m},
circular cross-section with radius \qty{0.005}{m},
elasticity modulus \qty{5.0}{MPa},
Poisson's ratio 0.5, and
mass density \qty{1100}{kg/m^3}.
The acceleration of gravity is assumed to be \qty{9.81}{m/s^2}.
The rod is discretized into 20 finite elements, yielding 21 nodes.
The time span simulated is \qty{2.48}{s}
and the time step is \qty{0.01}{s}.

\subsubsection{Oscillating beam}

This example considers a cantilever beam with neutral straight configuration
that is subject to the following simulation steps.
In the first step, the beam is statically deformed
under the action of a moment at the free end.
The concentrated load applied is such that the cantilever beam reaches
a deformed configuration corresponding to a semicircle.
In the second step, the cantilever beam is dynamically released
and oscillates exhibiting large displacements and large rotations
in the plane defined by the semicircle. There is no gravity.
The geometrical and material properties considered are the following:
length \qty{1.0}{m},
circular cross-section with radius \qty{0.005}{m},
elasticity modulus \qty{1.273}{GPa},
shear modulus \qty{0.637}{GPa}, and
mass density \qty{3183}{kg/m^3}.
The beam is uniformly discretized into 200 finite elements,
yielding 201 nodes.
The time span simulated is \qty{10}{s}
and the adopted time step is \qty{0.001}{s}.

\subsubsection{Triple pendulum}

This example considers a pendulum
that is built from three flexible straight links,
identical in their properties,
that are connected by means of spherical joints.
The triple pendulum is simply supported at one of its ends
and subject only to the action of its own weight induced by gravity.
Upon sudden release, it moves in a vertical plane
exhibiting nonlinear kinematics.
The geometrical and cross-sectional material properties
of each flexible link considered are the following:
length \qty{1.0}{m},
shear stiffness \Qty{5.0e4}{N/m},
axial stiffness \Qty{1.0e5}{N/m},
bending stiffness \Qty{6.25e-1}{N/m^2},
torsional stiffness \qty{6.239}{N/m^2},
mass density per unit length \Qty{2.5e-1}{kg/m},
inertia density per unit length \Qty{1.562e-6}{kgm}.
The acceleration of gravity is assumed to be \qty{9.81}{m/s^2}.
Each link is uniformly discretized into 20 finite elements,
yielding 21 nodes and thus, a total of 63 nodes.
The time span simulated is \qty{10}{s}
and the adopted time step is \qty{0.001}{s}.

\subsubsection{Wind turbine}

This example considers the NREL \qty{5}{MW} wind turbine,
which comprises a tower, a nacelle, a hub, and three blades.
The nacelle and hub are modeled as rigid bodies
while the tower and blades are modeled as beams.
The tower is clamped at the bottom and at the top,
the nacelle is rigidly fixed to the tower.
The blades are rigidly fixed to the hub,
which is connected to the nacelle with a hinge,
and thus only allowed to rotate about its symmetry axis.
A force is applied in the middle of one blade.
The force grows linearly between 0 and \qty{5}{s},
reaching a maximum of \Qty{4.0e5}{N}.
After that, the force drops instantaneously to \qty{0}{N}
and remains so for the rest of the simulation.
In contrast to the three prior examples,
the motion developed by the wind turbine is truly three dimensional.
Details on the geometrical and material properties
are to be found in \cite{jonkman2009definition}.
The tower is discretized into 10 elements
while each blade is discretized into 48.
The nacelle and the hub are represented by a single node each.
The whole model accounts then for a total of 160 nodes.
The time span simulated is \qty{10}{s}
and the adopted time step is \qty{0.0025}{s}.


\subsection{Examples on $S^2$}
\label{sec:examples-S2}

To carry out our investigations on $S^2$, we select three examples:
a tumbling cylinder,
a free-flying plate, and
a shell pendulum.
All three examples have been used as well
for error statistics in \cref{sec:accuracy}.
As already stated before, the principal geodesic analysis on $S^2$
is also carried out for individual directors
of the triple pendulum and wind turbine examples
corresponding to the set of examples on $SO(3)$.

\subsubsection{Tumbling cylinder}

This example considers a cylinder that is free to move in space,
which has been already investigated many times in the past,
see for instance \cite{Gebhardt2017} and references therein.
A complex set of forces is applied at different locations of the cylinder.
This set of forces is also subject to a time variation
that grows linearly in the first half of the initial second
and decreases linearly in its second half.
After that, the forces remain zero for the rest of the simulation.
The motion developed by the tumbling cylinder is truly three dimensional
and exhibits large displacements and rotations.
The geometrical and material properties of the cylinder are the following:
radius \qty{7.5}{m},
height \qty{3.0}{m},
thickness \qty{0.02}{m},
elastic modulus \qty{0.2}{GPa},
Poisson's ratio 0.25, and
mass density \qty{1.0}{kg/m^3}.
The tumbling cylinder is uniformly discretized into 48 elements
along its perimeter and into 3 along its height.
Thus it has 144 elements in total yielding 68 nodes,
which means that the mesh has a seam modeled as continuous connection.
The time span simulated is \qty{3.0}{s}
and the adopted time step is \qty{0.01}{s}.

\subsubsection{Free-flying plate}

This example considers a plate that is free to move in space,
which has been as well investigated before by Gebhardt and Rolfes
\cite{Gebhardt2017} among others.
A complex set forces is applied at different locations of the plate.
This set of forces is also subject to a time variation
that grows linearly between 0 and \qty{0.002}{s}
and decreases linearly between 0.002 and \qty{0.004}{s}.
After that, the forces remain zero for the rest of the simulation.
The motion of the free-flying plate shows large displacements and rotations.
Thus, its nonlinear kinematics is apparent.
Moreover, the dynamical behavior is very rich in high-frequency content.
The geometrical and material properties of the plate are the following:
length \qty{0.3}{m},
width \qty{0.06}{m},
thickness \qty{0.002}{m},
elastic modulus \qty{206.0}{GPa},
Poisson's ratio 0.0, and
mass density \qty{7300.0}{kg/m^3}.
The free-flying plate is uniformly discretized into 30 elements
along its length and into 4 along its width, \ie,
120 elements in total and yielding 155 nodes.
The time span simulated is \qty{0.1}{s}
and the adopted time step is \qty{0.0001}{s}.

\subsubsection{Shell pendulum}

This example considers a shell pendulum
that is simply supported at one of its edges
and subject only to the action of its own weight induced by gravity.
Upon sudden release,
the pendulum moves in space exhibiting nonlinear kinematics,
but with the transversal plane as symmetry plane for the motion.
The geometrical and material properties of the pendulum are the following:
length \qty{300.0}{m},
width \qty{60}{m},
thickness \qty{0.02}{m},
elastic modulus \qty{2.06}{MPa},
shear modulus \qty{1.03}{MPa}, and
mass density \qty{780}{kg/m^3}.
Some of the properties are intentionally chosen to be nonrealistic
to the purpose only of triggering a complex dynamical behavior.
The acceleration of gravity is assumed to be \qty{9.81}{kg/s^2}.
The shell pendulum is uniformly discretized into 30 elements
along its length and into 4 along its width, \ie,
120 elements in total and yielding 155 nodes.
The time span simulated is \qty{3.098}{s}
and the adopted time step is \qty{0.001}{s}.


\subsection{Computational Experiments}
\label{sec:experiments}

For each example described above the dynamic simulation
yields several hundred to several thousand snapshots with
several dozen to several hundred directors each;
see \cref{tab:examples}.
Each simulation is performed with a relative numerical accuracy
between \num{1e-10} and \num{1e-8}.
The configuration manifold is either
$(\R^3 \x \SO3)^{N_n}$ or $(\R^3 \x S^2)^{N_n}$
where $N_n$ is the number of nodes,
and we are only interested in the PGA
on $\SO3^{N_n}$ respectively $(S^2)^{N_n}$,
with lifts in $\so3^{N_n}$ respectively
$\prod_{i=1}^{N_n} (T_{d_{i0}}S^2)^{N_n}$.
\begin{table}
  \small
  \caption{Overview of examples; the snapshot matrix $Y$ has
    dimensions $m = \text{dim(lift)}$ and $n = \text{snapshots}$.}
  \label{tab:examples}
  \begin{tabular}{lrrrrr}
    \toprule
    Dynamic system & motion & nodes & directors & dim(lift) & snapshots \\
    \midrule
    Rubber rod        & 2d &  21 &  63 &  63 &   249 \\
    Oscillating beam  & 2d & 201 & 603 & 603 & 10001 \\
    Triple pendulum   & 2d &  63 & 189 & 189 & 10001 \\
    Wind turbine      & 3d & 160 & 480 & 480 &  4001 \\
    \midrule
    Tumbling cylinder & 3d &  68 &  68 & 116 &   301 \\
    Free-flying plate & 3d & 155 & 155 & 310 &  1001 \\
    Shell pendulum    & 3d (nearly 2d) & 155 & 155 & 310 &  3099 \\
    \bottomrule
  \end{tabular}
\end{table}

In order to evaluate the impact of lift-and-project errors
on physical quantities associated with the resulting trajectories,
we compute the total energy
\emph{i)} along the original sequence of snapshots;
\emph{ii)} along the lifted and projected
sequence of snapshots with matrix $Y$ \eqref{eq:lift:Y}; and
\emph{iii), iv)} along two sequences of snapshots
obtained by projecting truncated singular value decompositions $Y_p$
after the lifting \eqref{eq:trsvd}, \eqref{eq:tr-proj}:
a ``good'' one with a larger rank $p$, and
a ``poor'' one with a smaller rank $p$.
Relative differences of each lift-and-project energy
to the original total energy (which is invariant over time)
are presented in \cref{fig:energy:so3} for the three \SO3 examples
and in \cref{fig:energy:s2} for the two $S^2$ examples.
The values of the reference energy (denominators)
of the relative differences are as follows.
For the triple pendulum (\qty{10.9}{J})
and the shell pendulum (\qty{0.0406}{J}),
the reference energy is the difference between
initial and minimal potential energy reached during the simulation,
\ie, the maximal sum of kinetic energy plus elastic energy.
For the wind turbine (\Qty{2.62e5}{J}),
the tumbling cylinder (\qty{445}{J}),
and the free-flying plate (\qty{246}{J}),
it is the energy added to each system by applying the respective initial forces.

\newcommand\invimg[1]{\twtrimg[0.49]{15 20 46 27}{#1.png}}
\begin{figure}
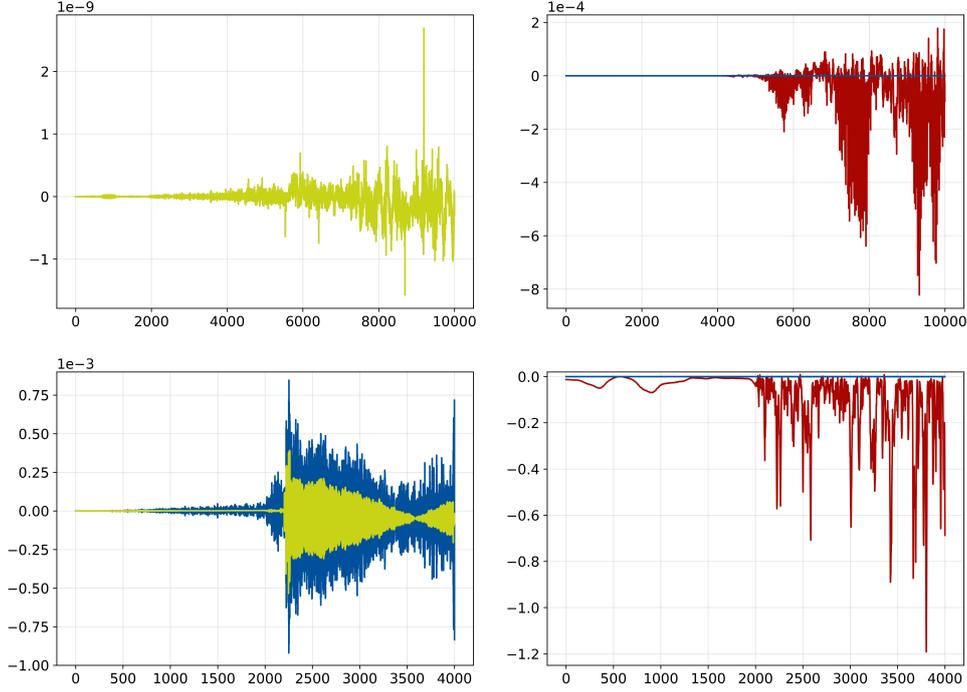

  \centering
  \invimg{pend_diff_orig_lift_good}\hfill
  \invimg{pend_diff_lift_good_bad}\\[2ex]
  \invimg{wind_diff_orig_lift_good}\hfill
  \invimg{wind_diff_lift_good_bad}
  \caption{\SO3: relative difference of original vs.\ \explift energy.
    Triple pendulum
    (top left: $Y$ green, $Y_{63}$ blue; right: $Y_{63}$ blue, $Y_{62}$ red),
    wind turbine
    (bottom left: $Y$ red, $Y_{345}$ blue; right: $Y_{345}$ blue, $Y_{35}$ red).}
  \label{fig:energy:so3}
\end{figure}

\begin{figure}
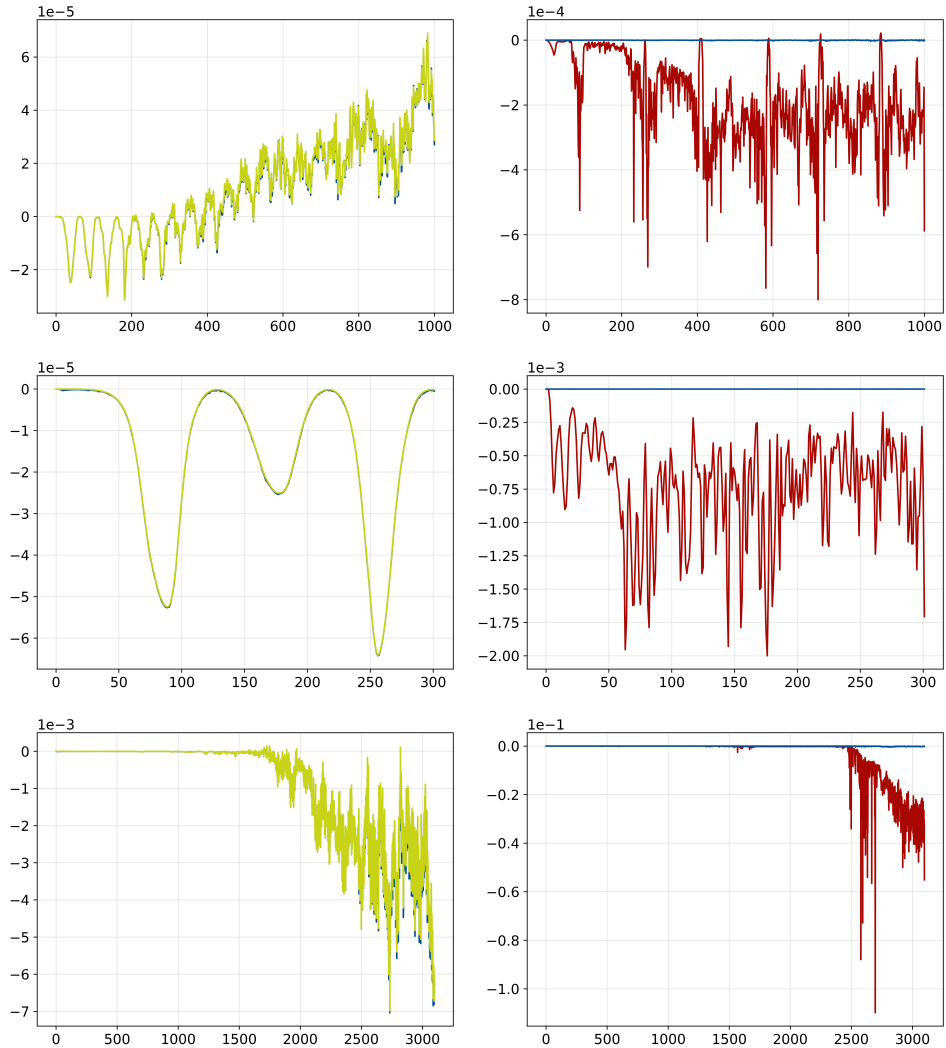

  \centering
  \invimg{ffp_diff_orig_lift_good}\hfill
  \invimg{ffp_diff_lift_good_bad}\\[2ex]
  \invimg{tc_diff_orig_lift_good}\hfill
  \invimg{tc_diff_lift_good_bad}\\[2ex]
  \invimg{sp_diff_orig_lift_good}\hfill
  \invimg{sp_diff_lift_good_bad}%
  \caption{$S^2$: relative difference of original vs.\ \explift energy.
    Tumbling cylinder (top left:
    $Y$ green, $Y_{84}$ blue; right: $Y_{84}$ blue, $Y_{47}$ red)
    free-flying plate (middle left:
    $Y$ green, $Y_{295}$ blue; right: $Y_{295}$ blue, $Y_{222}$ red),
    shell pendulum\ (bottom left:
    $Y$ green, $Y_{239}$ blue; right: $Y_{239}$ blue, $Y_{133}$ red).}
  \label{fig:energy:s2}
\end{figure}

For the triple pendulum we see in \cref{fig:energy:so3} (top left)
that the lift-and-project map without truncation (green)
and with truncation at $p = 63$ (blue, invisible)
give visually identical results,
with a maximal error of about \num{3e-8} in conserving the original energy.
This is because we have a planar rotation that is exactly reproduced
by the selected 63 (out of 189) SVD modes:
the relative truncation error is numerically zero at
$\sigma_{63} / \sigma_1 = 2 \epsM \approx \num{4e-16}$
and $\mnorm{Y - Y_{63}} / \sigma_1 \approx \num{8e-17}$.
With just one mode less, at $p = 62$, the energy error grows drastically
to a maximum of about \num{9e-3} (red curve in top right plot).
Here we have truncation errors
$\sigma_{62} / \sigma_1 \approx \num{1e-4}$ and
$\mnorm{Y - Y_{62}} / \sigma_1 \approx \num{2e-6}$.
As a function of time, all errors start out very small
and start growing significantly only after about 5000 time steps.

The 3d motion of the wind turbine exhibits significantly larger energy errors
(\cref{fig:energy:so3}, bottom left).
The maximal error without truncation is roughly {5e-4} (green curve),
and it doubles approximately
with a good approximation by 345 (out of 480) SVD modes (blue curve).
Here the truncation errors are
$\sigma_{345} / \sigma_1 \approx \num{1e-8}$ and
$\mnorm{Y - Y_{345}} / \sigma_1 \approx \num{2e-9}$.
The poor approximation at $p = 35$ has truncation errors
similar to the poor approximation of the triple pendulum,
$\sigma_{35} / \sigma_1 \approx \num{4e-5}$ and
$\mnorm{Y - Y_{35}} / \sigma_1 \approx \num{2e-6}$.
This gives very large energy errors up to 1.2 or 120\%
(red curve, bottom right).
The temporal behavior is similar to the triple pendulum,
with rather small errors during the first half of the simulation.

For the shell pendulum in \cref{fig:energy:s2} (bottom),
results are in part comparable to the triple pendulum:
the lift-and-project energy error without truncation (green)
and with a good aprroximation at $p = 239$ (blue)
are visually almost identical,
both with a maximum of roughly \num{7e-3}
and with very small values for the first half of the simulation.
Here the truncation error is $\sigma_{\new{239}} / \sigma_{\new1} \approx \num{7e-5}$.
The poor approximation, however, with only 133 (out of 310) modes
and a truncation error of $\sigma_{\new{133}} / \sigma_{\new1} \approx \num{1e-3}$,
yields comparatively large energy errors up to 0.11, \ie, 11\%.
On the other hand, these large errors develop later on
during the simulated time span.

The energy errors of the free-flying plate and of the tumbling cylinder
in \cref{fig:energy:s2} (top and middle) behave similarly,
except that in both cases
the errors start developing early on during the simulation
while accuracies of the lift and of the selected approximations are higher
than for the shell pendulum.
The truncation errors for the free-flying plate are
$\sigma_1 / \sigma_{295} \approx \num{1e-5}$ (good approximation) and
$\sigma_1 / \sigma_{222} \approx \num{1e-4}$ (poor approximation).
For the tumbling cylinder they are
$\sigma_1 / \sigma_{84} \approx \num{8e-7}$ (good approximation) and
$\sigma_1 / \sigma_{47} \approx \num{5e-5}$ (poor approximation).

\subsection{Discrete Trajectories and Lifts}
\label{eq:discrete-lifts}

Finally we illustrate some properties of the logarithms and lifts
for selected discrete trajectories of the examples.
Here we call the rotation matrices
of the finite element nodes \emph{orientations}
to distinguish them from the rotations that the physical bodies perform.

\begin{figure}
  \centering
  \twtrimg[0.49]{130 60 80 105}{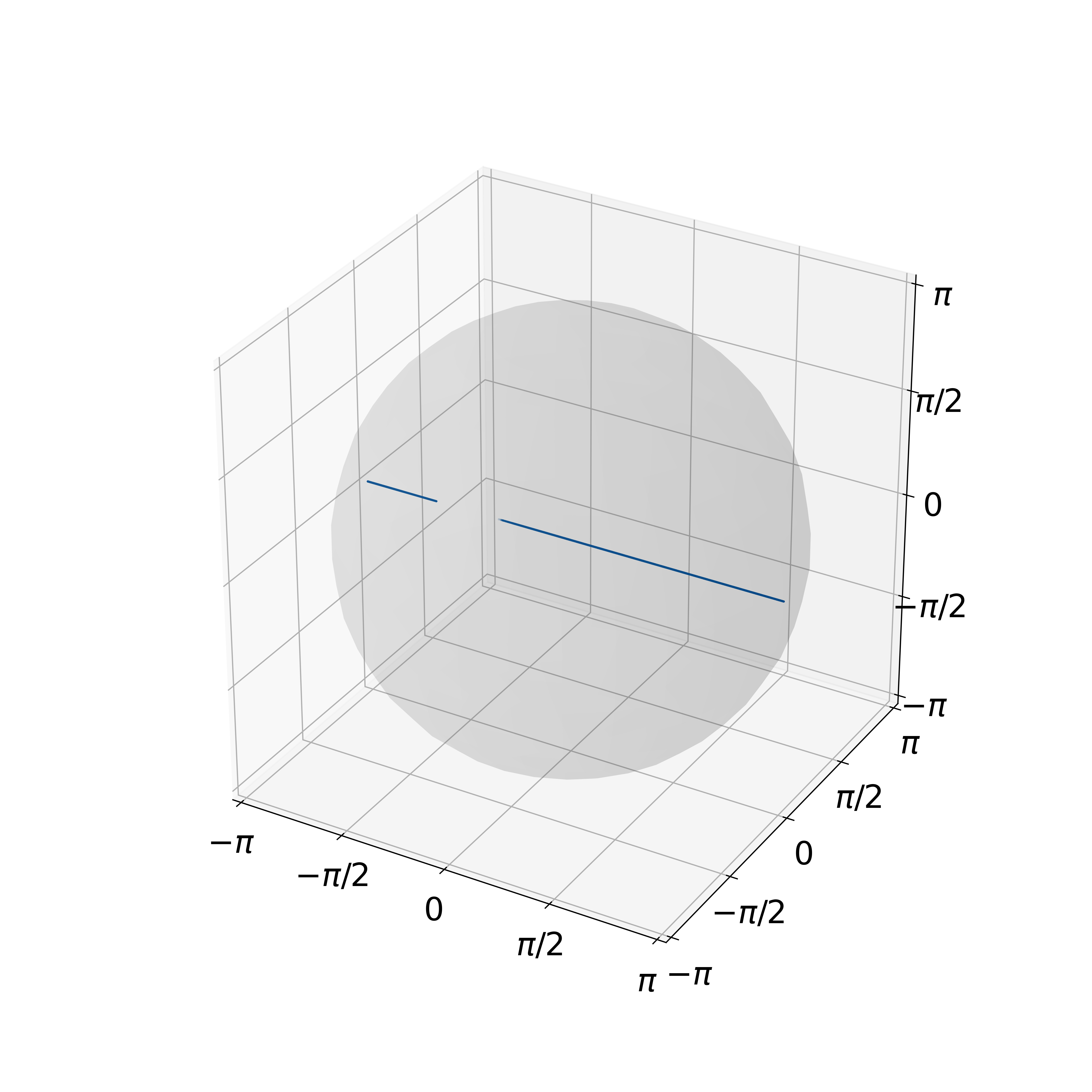}\hfill
  \twtrimg[0.49]{130 60 80 105}{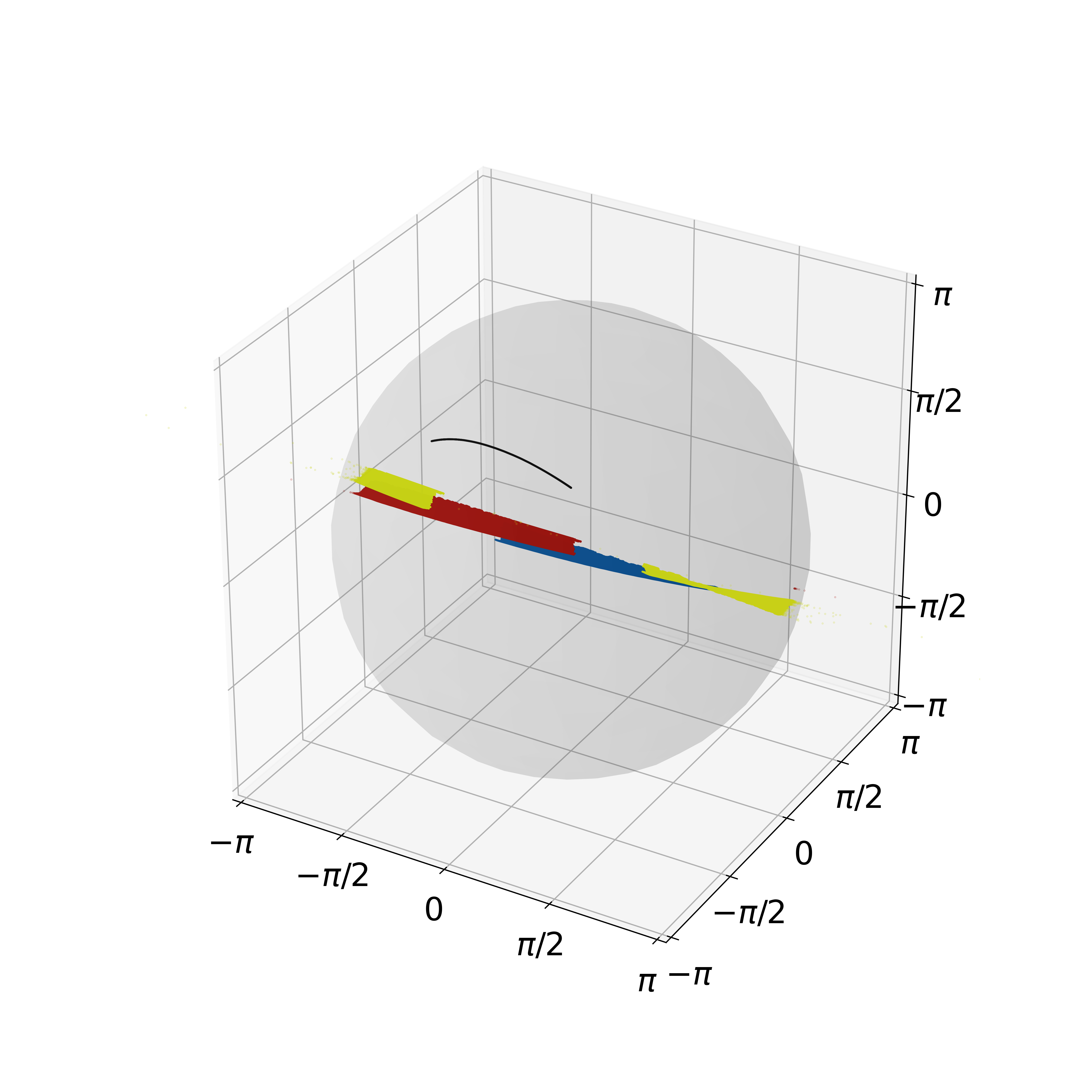}%
  \caption{Distribution of logarithms of rotations in Lie algebra \so3:
    triple pendulum (left), wind turbine (right: tower black, blades colored).}
  \label{fig:dist:so3-log}
\end{figure}

In \cref{fig:dist:so3-log} we see the logarithms
of all $\num{630063} = 63 \times 10001$
orientations of the triple pendulum (left)
and all $\num{640160} = 160 \times 4001$
orientations of the wind turbine (right),
which lie in the closed ball $\_B_\pi \subset \so3$ (indicated in grey).
Because of the planar motion of the triple pendulum,
all logarithms lie in a one-dimensional subspace of \so3,
\ie, on a straight line, which corresponds to the fixed
direction that the axes of rotation of all nodes share.
This is immediately apparent from the picture.
For the wind turbine, it is clearly visible that the logarithms
of all trajectories of blades' orientations (colored)
lie close to the same straight line but with notable differences:
the blades' trajectories exhibit slight deviations
from strictly planar rotation, so that a true 3d motion results.
The logarithms of the tower orientations (black)
show even larger deviations from a planar motion.

For the $S^2$ examples, \cref{fig:dist:s2-log} (left) presents
the distributions of all individual directors on the sphere:
$\num{20468} = 68 \times 301$ for the tumbling cylinder,
$\num{155155} = 155 \times 1001$ for the free-flying plate, and
$\num{480345} = 155 \times 3099$ for the shell pendulum.
Here we observe quite different behavior:
the tumbling cylinder has 68 clearly distinguishable trajectories
of tightly spaced directors,
the free-flying plate has 155 trajectories of widely spaced directors
that cover the entire sphere almost uniformly,
and the shell pendulum has again 155 trajectories
of widely spaced directors that are, however,
concentrated along a great circle of the sphere.

\newcommand\spaceimg[1]{\thtrimg[0.27]{140 90 110 110}{#1.png}}
\newcommand\planeimg[1]{\thtrimg[0.27]{25 40 50 55}{#1.png}}
\begin{figure}
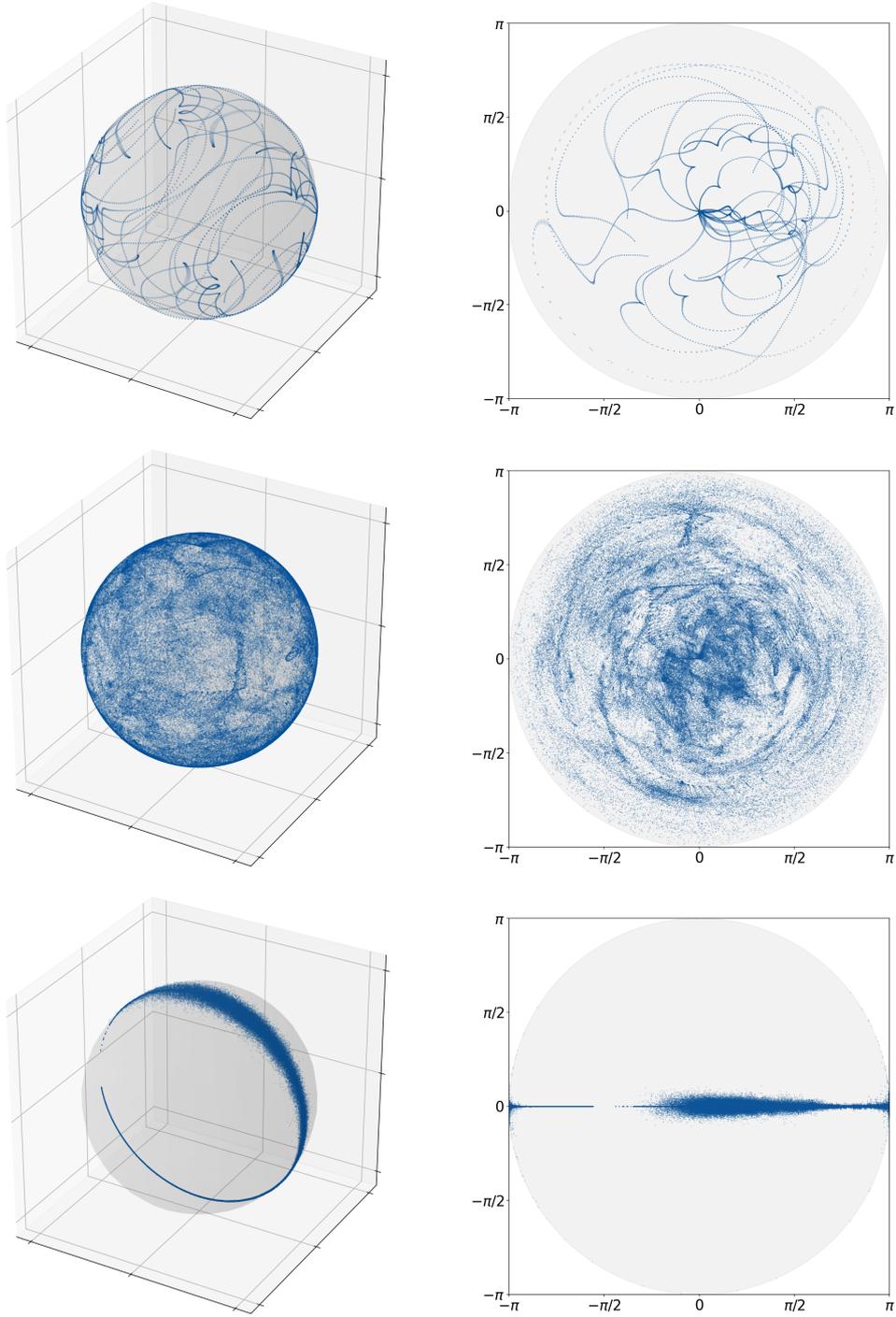

  \centering
  \spaceimg{tc_s2}\hfill\planeimg{tc_log_2d}\\[2ex]
  \spaceimg{ffp_s2}\hfill\planeimg{ffp_log_2d}\\[2ex]
  \spaceimg{sp_s2}\hfill\planeimg{sp_log_2d}
  \caption{Distribution of directors on sphere $S^2$ (left)
    and of their logarithms in tangent space $T_{\new{(0,0,1)}}S^2$ (right):
    tumbling cylinder (top), free-flying plate (middle),
    shell pendulum (bottom).}
  \label{fig:dist:s2-log}
\end{figure}

The corresponding logarithms in the closed disk
$\_D_\pi \subset T_{\new{(0,0,1)}}S^2$ (indicated in grey)
are displayed on the right of \cref{fig:dist:s2-log}.
Notice that the lifted discrete trajectory
of each finite element node $i \in \set{1, \dots, N_n}$
has its own tangent space, $T_dS^2$ with $d = d_{i0}$.
To illustrate all these lifted trajectories simultaneously,
we rotate each tangent space $T_dS^2$
into the tangent space $T_{(0,0,1)}S^2$ at the north pole via
\begin{align*}
  R(d) &\define \mat[3]{
    1 & 0 & -d_1\\
    0 & 1 & -d_2\\
    d_1 & d_2 & 1
  } \new{{}-{}} q(d) \mat[3]{
    d_1^2 & d_1d_2 & 0\\
    d_1d_2 & d_2^2 & 0\\
    0 & 0 & d_1^2+d_2^2
  }, &
  q(d) &\define \new{\frac1{1 + d_3}}.
\end{align*}
For a vector $v \in T_dS^2$
\new{we exploit the orthogonality relation $\sprod{d}{v} = 0$ to obtain}
\begin{equation*}
  R(d) v = \col(c){
    v_1 - d_1 v_3 \new{{}- q(d)} (v_1 d_1^2 + v_2 d_1 d_2) \\
    v_2 - d_2 v_3 \new{{}- q(d)} (v_1 d_1 d_2 + v_2 d_2^2) \\
    0
  } \new{{}= \col(c){v_1 - d_1 q(d) v_3 \\ v_2 - d_2 q(d) v_3 \\ 0}}.
\end{equation*}
(Here the coefficient map $q\: S^2 \new{{}\without{(0,0,-1)}} \to \R$
and consequently the rotation map $R\: S^2 \new{{}\without{(0,0,-1)}} \to \SO3$
\new{are undefined} at \new{the south pole} $d = (0,0,\new-1)$.)
Again the different behavior of the three dynamical systems
is immediately apparent from the pictures:
the logarithms of the directors of the tumbling cylinder
form clearly distinguishable trajectories,
the ones of the free-flying plate cover the disk $\_D_\pi$ almost uniformly,
and the logarithms of the shell pendulum
are concentrated along a straight line,
which indicates once more near-planar motion with slight deviations.

Analogous distributions for the individual directors $d_i$
of the wind turbine's orientations are shown in \cref{fig:dist:s2-wind},
again with all tangent spaces rotated into $T_{(0,0,1)}S^2$.
Here it is harder to infer the type of motion from the pictures.
However, except for the tower (black),
which clearly exhibits some true 3d motion,
the distributions of directors on the sphere (left)
clearly show that $d_1$ performs small movements around a nearly fixed director
for each node of each blade (colored);
this nearly fixed director is the horizontal axis of rotation.
Consequently, the directors $d_2$ and $d_3$ of all blades' nodes
move close to great circles.
Only the first logarithm plot, associated with~$d_1$,
appears to convey a clear situation.
However, due to the rotation of all tangent spaces into $T_{(0,0,1)}S^2$,
significant information gets lost
so that these pictures do not permit a well-founded interpretation.

\newcommand\logimg[1]{\thtrimg[0.27]{25 48 60 70}{#1.png}}
\newcommand\logbig[1]{\thtrimg[0.27]{10 10 35 35}
  {wind_neu_log_bigger_#1_2d.png}}
\begin{figure}
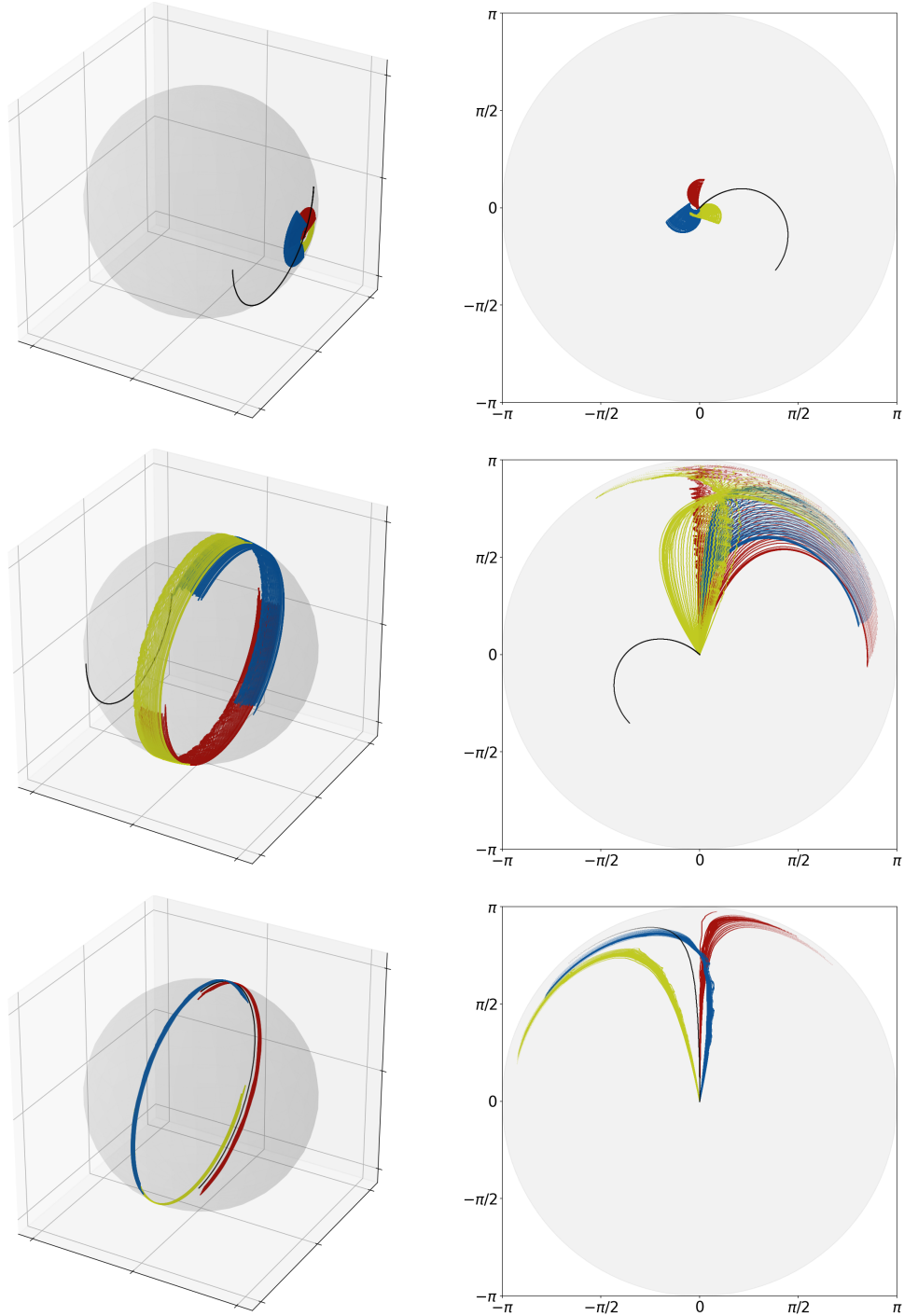

  \centering
  \spaceimg{wind_1_s2}\hfill\logimg{wind_neu_log_1_2d}\\[2ex]
  \spaceimg{wind_2_s2}\hfill\logimg{wind_neu_log_2_2d}\\[2ex]
  \spaceimg{wind_3_s2}\hfill\logimg{wind_neu_log_3_2d}%
  \caption{Rotations $R = [d_1 \ d_2 \ d_3] \in S^2 \x S^2 \x S^2$
    for wind turbine (tower black, blades colored):
    distribution of directors on sphere $S^2$ (left)
    and of their logarithms in tangent space $T_{(0,0,1)}S^2$;
    $d_1$ to $d_3$ (top to bottom).}
  \label{fig:dist:s2-wind}
\end{figure}

\begin{figure}
  \centering
  \twtrimg[0.49]{125 60 80 105}{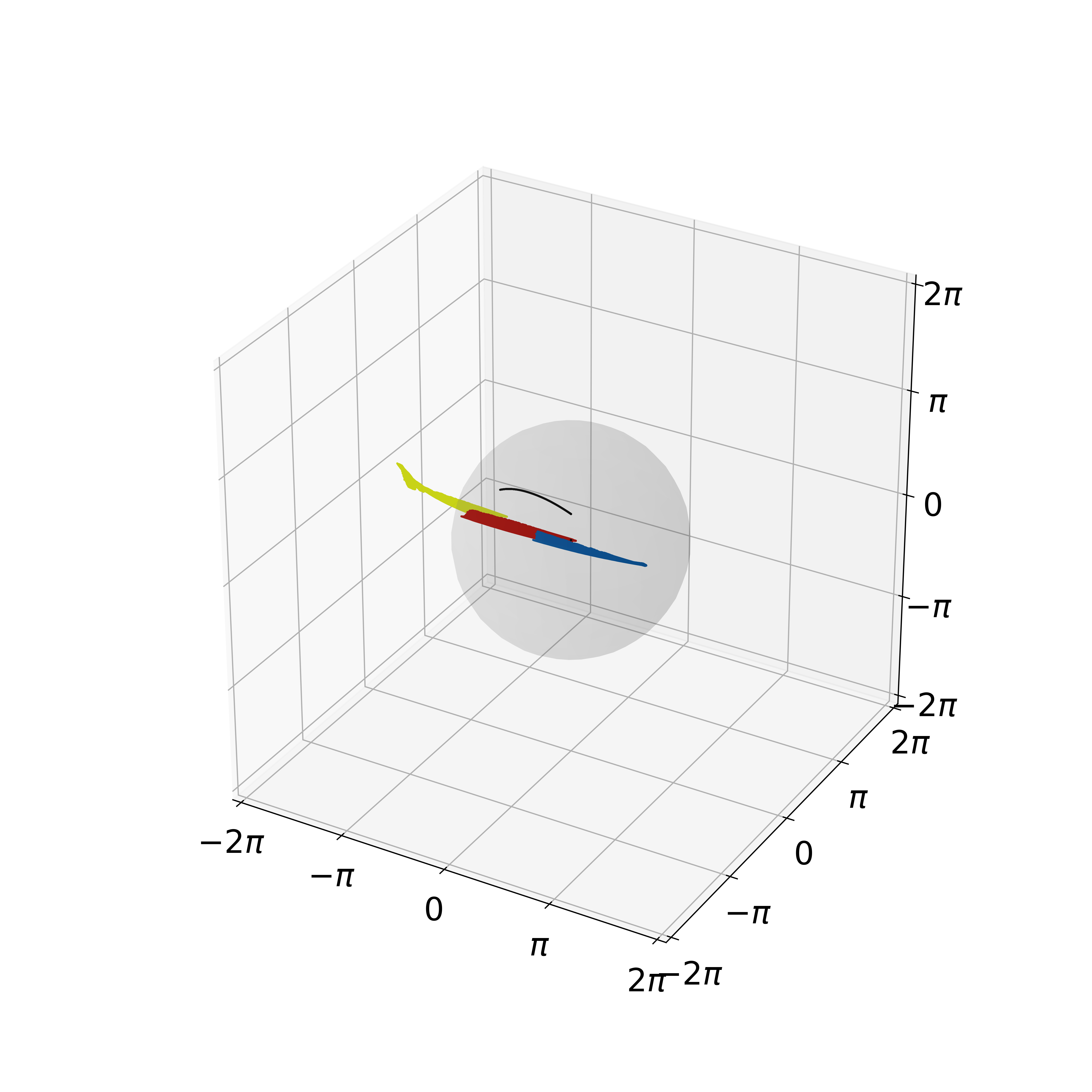}
  \caption{Lift with branch switching of wind turbine in \so3.}
  \label{fig:dist:lift-wind}
\end{figure}

Finally we present in \cref{fig:dist:lift-wind} a visualization
of the lifted trajectory of the wind turbine.
The ball $\_B_\pi$ is again indicated in grey
so that the presence of branch switchings for the logarithm map
is immediately apparent.
This illustration also shows that the lifts of all finite element nodes
of each blade cover a connected area,
which is again concentrated along a straight line.
In contrast, the green-yellowish blade appears to cover
two disconnected areas in \cref{fig:dist:so3-log}.

\section{Conclusions}
\label{sec:conclusion}

In this article, we presented new computational realizations of PGA
for the unit sphere $S^2$ and the special orthogonal group \SO3.
In addition, we carefully considered the consequences that singularities have
on the numerical accuracy of the proposed algorithms to construct
long-time smooth lifts across branches of the respective logarithmic maps.
Finally, we applied our PGA realization to investigate
several interesting mechanical systems.
Even provided that our results are encouraging,
especially regarding the potential for this approach
to be used in the derivation of structure-preserving reduced-order models,
there are still many open issues that will be considered in the sequel.

Among those open issues we identify as the most urgent ones:
\textit{i)} applying higher accuracy for lifting close to singularities,
which could be addressed by considering multiprecision arithmetic;
\textit{ii)} employing periodic, but discontinuous lifts
instead of continuous ones, but across branches;
\textit{iii)} carrying out PGA on the special Euclidean group $SE(3)$
instead of carrying out PCA on $\R^3$ and PGA on \SO3 separately.
Recall that the Lie groups $SE(3)$ and $\SO3 \times \R^3$ are diffeomorphic
as manifolds, but not isomorphic as groups;
\textit{iv)} further investigating degradation of mechanical invariants
such as momenta and energy upon truncation; and finally,
\textit{v)} deriving structure-preserving reduced-order models based on PGA.

\section*{Data availability}

The data supporting this manuscript can be provided upon reasonable request.

\section*{Conflict of interest statement}

The authors declare that they have no known competing financial interests or personal relationships that could have appeared to influence the work reported in this manuscript.

\section*{Acknowledgements}

The second and third authors are
funded by the Deutsche Forschungsgemeinschaft
(DFG, German Research Foundation) -- SFB1463 -- 434502799.

\bibliographystyle{siam}
\bibliography{lit}

\begin{thebibliography}{10}

\bibitem{Anderson_et_al:1992}
{\sc E.~Anderson, Z.~Bai, C.~Bischof, J.~Demmel, J.~Dongarra, J.~du~Croz,
  A.~Greenbaum, S.~Hammarling, A.~McKenney, S.~Ostrouchov, and D.~Sorensen},
  {\em LAPACK Users' Guide}, SIAM, Philadelphia, PA, 1992.
\newblock \url{http://www.netlib.org/lapack}.

\bibitem{Cazelles2018}
{\sc E.~Cazelles, V.~Seguy, J.~Bigot, M.~Cuturi, and N.~Papadakis}, {\em
  Geodesic {PCA} versus {Log-PCA} of histograms in the {Wasserstein} space},
  SIAM Journal on Scientific Computing, 40 (2018), pp.~B429--B456.

\bibitem{Chakraborty:2016}
{\sc R.~Chakraborty, D.~Seo, and B.~Vemuri}, {\em An efficient exact-{PGA}
  algorithm for constant curvature manifolds}, in Proceedings of the IEEE
  Conference on Computer Vision and Pattern Recognition, 2016, pp.~3976--3984.

\bibitem{Curry:2019}
{\sc C.~Curry, S.~Marsland, and R.~McLachlan}, {\em Principal symmetric space
  analysis}, Journal of Computational Dynamics, 6 (2019), pp.~251--276.

\bibitem{Fletcher:2003}
{\sc P.~Fletcher, C.~Lu, and S.~Joshi}, {\em Statistics of shape via principal
  geodesic analysis on {Lie} groups}, in 2003 IEEE Computer Society Conference
  on Computer Vision and Pattern Recognition, 2003. Proceedings, vol.~1, 2003,
  pp.~I--I.

\bibitem{Fletcher2004}
{\sc P.~Fletcher, C.~Lu, S.~M. Pizer, and S.~Joshi}, {\em Principal geodesic
  analysis for the study of nonlinear statistics of shape}, IEEE Transactions
  on Medical Imaging, 23 (2004), pp.~995--1005.

\bibitem{Fotouhi:2012}
{\sc H.~Fotouhi and M.~Golalizadeh}, {\em Exploring the variability of {DNA}
  molecules via principal geodesic analysis on the shape space}, Journal of
  Applied Statistics, 39 (2012), pp.~2199--2207.

\bibitem{Gebhardt2018cm}
{\sc C.~Gebhardt, B.~Hofmeister, C.~Hente, and R.~Rolfes}, {\em Nonlinear
  dynamics of slender structures: a new object-oriented framework},
  Computational Mechanics, - (2018), pp.~1--34.

\bibitem{Gebhardt2017}
{\sc C.~Gebhardt and R.~Rolfes}, {\em On the nonlinear dynamics of shell
  structures: combining a mixed finite element formulation and a robust
  integration scheme}, Thin-Walled Structures, 118 (2017), pp.~56--72.

\bibitem{Gebhardt:2020}
{\sc C.~Gebhardt, I.~Romero, and R.~Rolfes}, {\em {A new
  conservative/dissipative time integration scheme for nonlinear mechanical
  systems}}, Computational Mechanics, 65 (2020), pp.~405--427.

\bibitem{Gebhardt_et_al:2019}
{\sc C.~Gebhardt, M.~C. Steinbach, and R.~Rolfes}, {\em Understanding the
  nonlinear dynamics of beam structures: a principal geodesic analysis
  approach}, Thin-Walled Structures, 140 (2019), pp.~357--372.

\bibitem{Heard2006}
{\sc W.~B. Heard}, {\em Rigid Body Mechanics: Mathematics, Physics and
  Applications}, Wiley, 2006.

\bibitem{Heeren:2018}
{\sc B.~Heeren, C.~Zhang, M.~Rumpf, and W.~Smith}, {\em Principal geodesic
  analysis in the space of discrete shells}, Computer Graphics Forum, 37
  (2018), pp.~173--184.

\bibitem{Hernandez2021}
{\sc M.~Hernandez}, {\em Efficient momentum conservation constrained
  {PDE}-{LDDMM} with {G}auss–{N}ew{\-t}on–{K}ry{\-l}ov optimization,
  semi-{L}agrangian {R}unge–{K}utta solvers, and the band-limited
  parameterization}, Journal of Computational Science, 55 (2021), p.~101470.

\bibitem{Huckemann2010}
{\sc S.~Huckemann, T.~Hotz, and A.~Munk}, {\em Intrinsic shape analysis:
  geodesic {PCA} for {Riemannian} manifolds modulo isometric {Lie} group
  actions}, Statistica Sinica, 20 (2010), pp.~1--58.

\bibitem{Jolliffe:1986}
{\sc I.~Jolliffe}, {\em Principal Component Analysis}, Springer Series in
  Statistics, Springer New York, 1986.

\bibitem{jonkman2009definition}
{\sc J.~Jonkman, S.~Butterfield, W.~Musial, and G.~Scott}, {\em Definition of a
  5-{MW} reference wind turbine for offshore system development}, tech. rep.,
  National Renewable Energy Laboratory (NREL), Golden, CO., 2009.

\bibitem{Klingenberg:1982}
{\sc W.~Klingenberg}, {\em Riemannian Geometry}, vol.~1 of studies in
  mathematics, De Gruyter, Berlin, 1982.

\bibitem{Lazar:2017}
{\sc D.~Lazar and L.~Lin}, {\em Scale and curvature effects in principal
  geodesic analysis}, Journal of Multivariate Analysis, 153 (2017), pp.~64--82.

\bibitem{Lee:2007}
{\sc J.~Lee and M.~Verleysen}, {\em Nonlinear Dimensionality Reduction},
  Springer Series in Statistics, Springer New York, 2007.

\bibitem{Ren:2013}
{\sc P.~Ren, F.~Aziz, L.~Han, E.~Xu, R.~C. Wilson, and E.~R. Hancock}, {\em
  Geometricity and Embedding}, Springer London, 2013, pp.~121--155.

\bibitem{Said:2007}
{\sc S.~Said, N.~Courty, N.~Le~Bihan, and S.~Sangwine}, {\em Exact principal
  geodesic analysis for data on {SO(3)}}, in 2007 15th European Signal
  Processing Conference, 2007, pp.~1701--1705.

\bibitem{Salehian:2014}
{\sc H.~Salehian, D.~Vaillancourt, and B.~Vemuri}, {\em {iPGA}: incremental
  principal geodesic analysis with applications to movement disorder
  classification}, in International Conference on Medical Image Computing and
  Computer-Assisted Intervention, vol.~17, 2014, pp.~765--772.

\bibitem{Sittel:2017}
{\sc F.~Sittel, T.~Filk, and G.~Stock}, {\em Principal component analysis on a
  torus: theory and application to protein dynamics}, The Journal of Chemical
  Physics, 147 (2017), p.~244101.

\bibitem{Sommer:2014}
{\sc S.~Sommer, F.~Lauze, and M.~Nielsen}, {\em Optimization over geodesics for
  exact principal geodesic analysis}, Advances in Computational Mathematics, 40
  (2014), pp.~283--313.

\bibitem{Tournier_et_al:2009}
{\sc M.~Tournier, X.~Wu, N.~Courty, E.~Arnaud, and L.~Revéret}, {\em Motion
  compression using principal geodesics analysis}, Computer Graphics Forum, 28
  (2009), pp.~355--364.

\bibitem{Boost}
{\em {Boost C++ Libraries}}.
\newblock \url{http://www.boost.org/}, 1998-\the\year.

\end{thebibliography}

\end{document}